\documentclass[12pt]{article}
\title{Birational geometry of symplectic resolutions 
of nilpotent orbits}
\author{Yoshinori Namikawa}
\date{ }
\usepackage{amsmath,amssymb,amsthm, amscd}
\chardef\bslash=`\\
\newtheorem{Thm}{Theorem}[section]
\newtheorem{Cor}[Thm]{Corollary}
\newtheorem{Lem}[Thm]{Lemma}
\newtheorem{Prop}[Thm]{Proposition}
\newtheorem{Def}{Definition}
\newtheorem{Rque}[Thm]{Remark}
\newtheorem{Conj}{Conjecture}
\newtheorem{Exam}[Thm]{Example}

\def\cit{{\mathbb C}}

\def\0{{\mathcal O}}
\def\g{{\mathfrak g}}
\def\p{{\mathfrak p}}

\begin{document}
\maketitle

\section{Introduction}

Let $G$ be a complex simple Lie group and let 
$\g$ be its Lie algebra. Then $G$ has the adjoint 
action on $\g$. The orbit $\mathcal{O}_x$ of a 
nilpotent element $x \in \g$ is called a nilpotent 
orbit. A nilpotent orbit $\mathcal{O}_x$ admits 
a non-degenerate closed 2-form $\omega$ called 
the Kostant-Kirillov symplectic form. 
The closure $\bar{\mathcal{O}}_x$ of $\mathcal{O}_x$ 
then becomes a symplectic singularity. 
In other words, the 2-form $\omega$ extends to a 
holomorphic 2-form on a resolution of 
$\bar{\mathcal{O}}_x$. A resolution of 
$\bar{\mathcal{O}}_x$ is called a symplectic resolution 
if this extended form is everywhere non-degenerate 
on the resolution. A typical symplectic resolution 
of $\bar{\mathcal{O}}_x$ is obtained as the 
{\em Springer resolution} 
$$ T^*(G/P) \to \bar{\mathcal{O}}_x $$ 
for a suitable parabolic subgroup $P \subset 
G$. Here $T^*(G/P)$ is the cotangent bundle 
of the homogenous space $G/P$. 
Spaltenstein \cite{S} and Hesselink \cite{He} obtained a 
necessary and sufficient condition for 
$\bar{\mathcal{O}}_x$ to have a Springer resolution 
when $\g$ is a classical simple Lie algebra. 
Moreover, \cite{He} gave an explicit number of such 
parabolics $P$ up to conjugacy class that give 
Springer resolutions of $\bar{\mathcal{O}}_x$
(cf. \S 2). Recently, Fu \cite{Fu 1} has shown that 
every symplectic (projective) resolution 
is obtained as a Springer resolution. 
The following is one of main results of 
this paper. \vspace{0.12cm}

{\bf Theorem 4.4}. {\em Let 
$\mathcal{O}_x$ be a nilpotent orbit 
of a classical complex simple Lie algebra. 
Let $Y$ and $Y'$ be any two Springer resolutions 
of the closure $\bar{\mathcal{O}}_x$ of the 
nilpotent orbit. Then the birational map 
$Y --\to Y'$ can be decomposed into 
finite number of diagrams 
$Y_i \rightarrow X_i \leftarrow Y_{i+1}$ 
$(i = 1, ..., m-1)$ with $Y_1 = Y$ and 
$Y_m = Y'$ in such a way that each diagram 
is locally a trivial family of Mukai flops 
of type A or of type D.}   
\vspace{0.12cm}

A {\em Mukai flop of type A} is a kind of Springer 
resolutions; let $x \in \mathfrak{sl}(n)$ be 
a nilpotent element of type $[2^k, 1^{n-2k}]$ 
with $2k < n$. Then a Mukai flop of type A is 
the diagram of two Springer resolutions of 
$\bar{\mathcal{O}}_x$: 
$$ T^*G(k,n) \rightarrow \bar{\mathcal{O}}_x 
\leftarrow T^*G(n-k,n) $$ 
where $G(k,n)$ (resp. $G(n-k,n)$) is the 
Grassmannian which parametrizes $k$-dimensional 
(resp. $n-k$-dimensional) subspaces of 
$\mathbf{C}^n$. This flop naturally appears 
in the wall-crossing of the moduli spaces 
of various objects (eg. stable sheaves on 
K3 surfaces, quiver varieties and so on). 
On the other hand, a {\em Mukai flop of type D} 
comes from an orbit of a simple Lie algebra 
of type D. Let $x \in \mathfrak{so}(2k)$ be 
a nilpotent element of type $[2^{k-1}, 1^2]$, 
where $k$ is an odd integer with $k \geq 3$. 
Then $\bar{\mathcal{O}}_x$ admits two Springer 
resolutions 
$$ T^*G^+_{iso}(k,2k) \rightarrow \bar{\mathcal{O}}_x 
\leftarrow T^*G^-_{iso}(k,2k) $$ 
where $G^+_{iso}(k,2k)$ and $G^-_{iso}(k,2k)$ are 
two connected components of the orthogonal 
Grassmannian $G_{iso}(k,2k)$. For details on 
these Mukai flops, see \S 3. 
The factorizations in Theorem 4.4 make it 
possible to describe the ample cones and movable 
cones of symplectic resolutions of 
$\bar{\mathcal{O}}_x$. Moreover, Theorem 4.4 would 
clarify the geometric meaning of the results of 
Spaltenstein and Hesselink. To illustrate these, 
three examples will be given (see Examples 4.6, 
4.7, 4.8). 

Another purpose of this paper is to give an 
affimative answer to the following conjecture 
in the case of (the normalization of) a 
nilpotent orbit closure in a classical simple 
Lie algebra (Theorem 5.9).  \vspace{0.12cm}

{\bf Conjecture}(\cite{F-N}): {\em Let $W$ be 
a normal symplectic singularity. Then for any 
two symplectic resolutions $f_i: X_i \to 
W$, $i = 1,2$, there are deformations 
$\mathcal{X}_i \stackrel{F_i}\to \mathcal{W}$ 
of $f_i$ 
over a parameter space $S$  such that, for 
$s \in S - \{0\}$, $F_{i,s}: \mathcal{X}_{i,s} 
\to \mathcal{W}_s$ are isomorphisms. 
In particular, $X_1$ and $X_2$ are deformation 
equivalent. }   
\vspace{0.12cm}

This conjecture is already proved in \cite{F-N} when 
$W$ is a nilpotent orbit closure in 
$\mathfrak{sl}(n)$. On the other hand, a weaker 
version of this conjecture is proved in \cite{Fu 2} 
when $W$ is the normalization of a nilpotent 
orbit closure in a classical simple Lie algebra. 
According to the idea of Borho and Kraft \cite{B-K}, 
we shall define a deformation of 
$\bar{\mathcal{O}}_x$ by using a {\em Dixmier sheet}. 
Corresonding to each parabolic subgroup $P$, 
this deformation has a simultaneous resolution. 
These simultaneous resolutions would give the 
desired deformations of the conjecture. 
Details on the construction of them can be found 
in \S 5. \vspace{0.15cm}      

{\bf Notation and Convention}.  
A partition $\mathbf{d}$ of $n$ is a 
set of positive integers $[d_1, ..., d_k]$ 
such that $\Sigma d_i = n$ and $d_1 \geq d_2 
\geq ... \geq d_k$. We mean by $[d_1^{j_1}, 
..., d_k^{j_k}]$ the partition where $d_i$ 
appear in $j_i$ multiplicity.  
If $(p_1, ..., p_s)$ is a sequence of positive 
integers, then we define the partition 
$\mathbf{d} = \mathrm{ord}(p_1, ..., p_s)$ 
by $d_i := \sharp \, \{j; p_j \geq i\}$. 
In particular, for a partition $\mathbf{d}$, 
${ }^t\mathbf{d} := \mathrm{ord}
(d_1, ..., d_k)$ is called the dual partition 
of $\mathbf{d}$. We define $d^i := ({ }^td)_i$.

\section{Nilpotent orbits and  
Polarizations}

Let $G$ be a complex simple Lie group and let 
$\g$ be its Lie algebra. $G$ has the adjoint 
action on $\g$. The orbit $\mathcal{O}_x$ of a nilpotent 
element $x \in \g$ for this action is called 
a nilpotent orbit. This orbit carries a natural 
closed non-degenerate 2-form (Kostant-Kirillov 
form) $\omega$ (cf. \cite{C-G}, Prop. 1.1.5, \cite{C-M}, 1.3), 
and its closure 
$\bar{\mathcal{O}}_x$ becomes a {\em symplectic} singularity, 
that is, the symplectic 2-form $\omega$ extends to 
a holomorphic 2-form on a resolution $Y$ of 
$\bar{\mathcal{O}}_x$. When $\g$ is classical, $\g$ is naturally 
a Lie subalgebra of $\mathrm{End}(V)$ for a complex vector space 
$V$. Then we can attach a partition ${\bf d}$ of 
$n := \dim V$ to each orbit as the Jordan 
type of an element contained in the orbit. 
Here a 
partition ${\bf d} := [d_1, d_2,..., d_k]$ of $n$ is a set of 
positive integers with $\Sigma d_i = n$ and 
$d_1 \geq d_2  \geq ... \geq d_k$. When a  
number $e$ appears in the partition $\bf{d}$, we say 
that $e$ is a {\em part} of $\bf{d}$. We call $\bf{d}$ 
{\em very even} when $\bf{d}$ consists with only even 
parts, each having even multiplicity. 
The following result can be found, for example, in 
[C-M, \S 5].     

\begin{Prop} \label{nilprop}  
Let $\mathcal{N}(\g)$ be the set of nilpotent 
orbits of $\g$. 
\vspace{0.12cm}

(1)($A_{n-1}$): When $\g = \mathfrak{sl}(n)$, there is a 
bijection between $\mathcal{N}(\g)$ and 
the set of partitions $\bf{d}$ of $n$.   
\vspace{0.12cm} 
 
(2)($B_n$): When $\g = \mathfrak{so}(2n+1)$, there is a 
bijection between $\mathcal{N}(\g)$ and 
the set of partitions $\bf{d}$ of $2n+1$ such that 
even parts occur with even multiplicity.   
\vspace{0.12cm}

(3)($C_n$): When $\g = \mathfrak{sp}(2n)$, there is a 
bijection between $\mathcal{N}(\g)$ and the 
set of partitions $\bf{d}$ of $2n$ such that odd 
parts occur with even multiplicity
\vspace{0.12cm}

(4)($D_n$): When $\g = \mathfrak{so}(2n)$, there is a 
surjection $f$ from $\mathcal{N}(\g)$ to the set 
of partitions $\bf{d}$ of $2n$ such that even parts occur 
with even multiplicity. For a partition $\bf{d}$ which is 
not very even, $f^{-1}(\bf{d})$ consists of exactly one orbit, 
but, for very even $\bf{d}$, $f^{-1}(\bf{d})$ consists of exactly 
two different orbits.   
\end{Prop}  

Let $V$ be a finite dimensional ${\bf C}$-vector space.   
Then a flag $F := \{F_i\}_{1 \leq i \leq s}$ is a sequence 
of vector subspaces of $V$: 
$$  F_1 \subset F_2 \subset ... \subset F_s = V. $$ 
We put $p_i := \dim (F_i/F_{i-1})$ and the sequence 
$(p_1, ..., p_s)$ is called the type of $F$. 
The following is well-known (cf. \cite{He}, Lemma 4.3):  

\begin{Prop} \label{para A} 
Every parabolic subgroup $P$ of $SL(V)$ is 
a stabilizer of a suitable flag $F$ of $V$, 
and conversely, for every flag 
$F$ of $V$, its stabilizer group $P \subset 
SL(V)$ is a parabolic subgroup. The conjugacy 
class of a given parabolic subgroup is completely 
determined by the type of the corresponding flag.   
\end{Prop}

Let $\epsilon$ denote the number $0$ or $1$. 
Assume that $V$ is a ${\bf C}$-vector space equipped 
with a non-degenerate bilinear form $<,>$ such that 
$$ <v, w> = (-1)^{\epsilon} <w, v>,  (v,w \in V). $$  
When $\epsilon = 0$ (resp. $\epsilon = 1$), 
this means that the bilinear form 
is symmetric (resp. skew-symmetric). 
We put

$$ H := \{x \in GL(V); <xv, xw> = <v, w>, (v,w \in 
V)\}, $$ 
and 
$$ G := \{x \in H; \mathrm{det}(x) = 1\}. $$ 

Note that 
$$H = \left\{ 
\begin{array}{rl}    
O(V) & \quad (\epsilon = 0)\\ 
Sp(V) & \quad (\epsilon = 1) 
\end{array}\right.$$ and 

$$G = \left\{ 
\begin{array}{rl}    
SO(V) & \quad (\epsilon = 0)\\ 
Sp(V) & \quad (\epsilon = 1) 
\end{array}\right.$$

A flag $F := \{F_i\}_{1 \leq i 
\leq s}$ of $V$ is called {\em isotropic} if 
$F_i^{\perp} = F_{s-i}$ for $1 \leq i \leq s$. 
An isotropic flag $F$ is {\em admissible} if the 
stabilizer group $P$ of $F$ has no finner 
stabilized flag than $F$. In other words, let 
$P_F := \{g \in G; gF_i \subset F_i \: \forall i\}$. 
Then, for any $i$, there is no $P_F$-invariant 
subspace $F'_i$ 
such that $F_i \subset F'_i \subset F_{i+1}$ with 
$F'_i \neq F_i, F_{i+1}$. 
When the length $s$ of an isotropic flag $F$ is even, one can 
write the type of $F$ as $(p_1, ..., p_k, p_k, ..., p_1)$ 
with $k = s/2$. On the other hand, when $s$ is odd, 
one can write the type of $F$ as $(p_1, ..., p_k,q, p_k, 
..., p_1)$ with $k = (s-1)/2$. 

\begin{Exam}      
Assume that $\epsilon = 0$ and  
$q = 2$. Then one can always find 
a $P_F$-invariant subspace $F'_k$ such that 
$F_k \subset F'_k \subset F_{k+1}$ and $\dim (F'_k/F_k) 
= 1$. Therefore, $F$ is an isotropic flag which is not 
admissible. This is the only case where an isotropic flag 
is not admissible. 
\end{Exam}

\begin{Prop} \label{para B,C,D}
Let $V$ be a finite dimensional $\bf C$-vector 
space with a non-degenerate bilinear form 
and let $\epsilon$ be $0$ or $1$ according as 
the bilinear form is symmetric or skew-symmetric. 
Let $G$ and $H$ be the same as above.
\vspace{0.12cm}

(i) Any parabolic subgroup $P$ of $G$ 
is a stabilizer group of an admissible flag 
of $V$, and conversely, the stabilizer group 
of an admissible flag becomes a parabolic 
subgroup of $G$. 
\vspace{0.12cm}

(ii) There is a one-to-one correspondence 
between $H$-conjugacy classes of parabolic 
subgroups of $G$, and types of the stabilized 
admissible flags.  
\vspace{0.12cm}

(iii) When $\epsilon = 1$, 
the $H$-conjugacy class of a parabolic subgroup of $G$ 
coincides with the $G$-conjugacy class of a parabolic 
subgroup of $G$. 
\vspace{0.12cm}

(iv) When $\epsilon = 0$, 
the $H$-conjugacy class of a parabolic subgroup of $G$ 
coincides with the $G$-conjugacy class of a parabolic 
subgroup of $G$ except in the case where the type of 
the stabilized flag satisfies $q = 0$ and $p_k \geq 2$
(cf. the notation above). 
In this particular case, the $H$-conjugacy class splits 
into two $G$-conjugacy classes.   
\end{Prop}  

\begin{Def}
If a parabolic subgroup of $G$ has 
a stabilized (admissible) flag $F$ of type  
$(p_1, ..., p_k,q,p_k, ..., p_1)$, then 
$\pi := \mathrm{ord}(p_1, ..., p_k, q, p_k, ..., p_1)$ 
is called the {\em Levi type} of $P$. 
\end{Def}

\begin{Rque}
Assume that $\epsilon = 0$  
and $\dim V$ is even. Put $p = 1/2\dim V$ and 
let $G_{\mathrm{iso}}(p,V)$ be the {\em orthogonal 
Grassmannian} which parametrizes isotropic $p$ dimensional 
subspaces of $V$. The orthogonal Grassmannian has two 
connected components. Take two points $[F_1]$, $[F_2]$ 
from different components of $G_{\mathrm{iso}}(p,V)$. 
Then one can find $h \in H$ such that $h[F_1] = [F_2]$, 
but, when $p \geq 2$, there are no such elements in 
$G$. This is the reason why the exceptional cases occur 
in (iv) of the proposition   
\end{Rque}

Let $G$ be a complex simple Lie algebra as above and let 
$x \in \g$ be a nilpotent element. Then a 
parabolic subgroup $P$ of $G$ is called a 
{\em polarization} of $x$ if $x \in \mathfrak{n}(P)$ 
and $\dim {\mathcal O}_x = 2 \dim G/P$, 
where $\mathfrak{n}(P)$ is 
the nil-radical of $\p := \mathrm{Lie}(P)$. 
If we have a polarization $P$ of $x$, then we can 
define a map 
$$\mu: G \times_P \mathfrak{n}(P) \to \bar{\mathcal{O}}_x$$ 
as $\mu ([g,x]) = \mathrm{Ad}_g(x)$. Here $G \times_P 
\mathfrak{n}(P)$ is the quotient space of 
$G \times \mathfrak{n}(P)$ by the equivalence relation:  
$$ (g,x) \sim (g', x') \Leftrightarrow g' = gp, x' = 
\mathrm{Ad}_{p^{-1}}(x), \exists p \in P. $$ 
$G \times_P 
\mathfrak{n}(P)$ is a vector bundle over $G/P$ and 
it coincides with the cotangent bundle $T^*(G/P)$ 
of $G/P$. The map $\mu$ is a generically finite, proper 
and surjective map. If $\mathrm{deg}(\mu) = 1$, then 
$\mu$ becomes a resolution of $\bar{\mathcal{O}}_x$. 
In this case, $\mu$ is called the {\em Springer resolution} 
of $\bar{\mathcal{O}}_x$ with respect to $P$.   
\vspace{0.12cm}

Let $x \in \g$ be a nilpotent orbit and denote by 
$\mathrm{Pol}(x)$ the set of polarizations of $x$.  

\begin{Thm} \label{Pol A}   
Let $x \in \mathfrak{sl}(n)$ be a nilpotent 
element. Then $\mathrm{Pol}(x) \neq \emptyset.$ 
Assume that $x$ is of type 
$\mathbf{d} = [d_1, ..., d_k]$. Then 
$P \in \mathrm{Pol}(x)$ has the flag type 
$(p_1, ..., p_s)$ such that $\mathrm{ord}(p_1, ..., p_s) 
= \mathbf{d}$. Conversely, for any sequence 
$(p_1, ..., p_s)$ with 
$\mathrm{ord}(p_1, ..., p_s) 
= \mathbf{d}$, there is a unique polarization 
$P \in \mathrm{Pol}(x)$ which has the flag type 
$(p_1, ..., p_s)$. 
\end{Thm}

{\em Proof}. We shall construct a flag $F$ of type $(p_1, ..., p_s)$ 
such that $xF_i \subset F_{i-1}$ for all $i$. 
We identify the partition {\bf d} with a Young table  
consisting of $n$ boxes, where the $i$-th row consists of $d_i$ boxes 
for each $i$. We denote by $(i,j)$ the box of {\bf d} lying on the 
$i$-th row and on the $j$-th column. 
Let $e(i,j)$, $(i,j) \in {\mathbf d}$ be a {\em Jordan basis} of $V := 
\mathbf{C}^n$ such that $xe(i,j) = e(i-1,j)$. 
We consruct a flag by the induction on $n$. 
Define first $F_1 := \Sigma_{1 \leq j \leq p_1} \mathbf{C} e(1,j)$. 
Then $x$ induces a nilpotent endomorphism $\bar{x}$ of 
$V/F_1$. The Jordan type of $\bar{x}$  
is $[d_1-1, ..., d_{p_1}-1, d_{p_1+1}, ..., d_k]$. 
Note that this coincides with $\mathrm{ord}(p_2, ..., p_k)$. 
By the induction hypothesis, we already have a flag of type 
$(p_2, ..., p_k)$ on 
$V/F_1$ stabilized by $\bar{x}$; hence we have a desired 
flag $F$.  Let $P$ be the stabilizer group of $F$. Then 
it is clear that $x \in \mathfrak{n}(P)$. By an explicit calculation 
$\dim \mathcal{O}_x = 2 \dim G/P$.    
\vspace{0.12cm}

Next consider simple Lie algebras of type 
$B$, $C$ or $D$. Let $V$ be an $n$ dimensional 
${\bf C}$-vector space with a non-degenerate 
symmetric (skew-symmetric) form. As above, $\epsilon 
= 0$ when this form is symmetric and $\epsilon = 1$ 
when this form is skew-symmetric.  
Let $P_{\epsilon}(n)$ be the set of 
partitions $\mathbf{d}$ of $n$ such that 
$\sharp\{i; d_i = m\}$ is even for every integer 
$m$ with $m \equiv \epsilon$ (mod 2). Note that 
these partitions are nothing but those which appear 
as the Jordan types of nilpotent elements of 
$\mathfrak{so}(n)$ or of $\mathfrak{sp}(n)$.  
Next, let $q$ be a non-negative integer and assume 
moreover that $q \ne 2$ when $\epsilon = 0$. 
We define $\mathrm{Pai}(n,q)$ to be the set of 
partitions $\mathfrak{\pi}$ of $n$ such that 
$\pi_i \equiv 1$ (mod 2) if $i \leq q$ and 
$\pi_i \equiv 0$ (mod 2) if $i > q$. Note 
that, if $(p_1, ..., p_k,q, p_k, ..., p_1)$ 
is the type of an admissible flag of $V$, then 
$\mathrm{ord}(p_1, ..., p_k,q, p_k, ..., p_1) 
\in \mathrm{Pai}(n,q)$. 
Now we shall define the {\em Spaltenstein map} $S$  
from $\mathrm{Pai}(n,q)$ to $P_{\epsilon}(n)$. 
For $\pi \in \mathrm{Pai}(n,q)$, let 
$$ I(\pi) := \{ j \in \mathbf{N}\vert j \not\equiv n \,(\mathrm{mod}\,2), 
\pi_j \equiv \epsilon \,(\mathrm{mod}\,2), 
\pi_j \geq \pi_{j+1} + 2\}. $$ 
Then the Spaltenstein map     
$$ S: \mathrm{Pai}(n,q) \to P_{\epsilon}(n)$$ 
is defined as 
$$S(\pi)_j := \left\{ 
\begin{array}{rl}    
\pi_j -1 & \quad (j \in I(\pi))\\ 
\pi_j +1 & \quad (j-1 \in I(\pi)) \\
\pi_j & \quad (otherwise)
\end{array}\right.$$       

\begin{Thm}
Let $G$ be $SO(V)$ or $Sp(V)$ according as $\epsilon = 0$ 
or $\epsilon = 1$. Let $x \in \g$ be a nilpotent element 
of type $\mathbf{d} \in P_{\epsilon}(n)$. For $\pi \in 
\mathrm{Pai}(n,q)$ , define $\mathrm{Pol}(x, \pi)$ to be 
the set of polarizations of $x$ with Levi type $\pi$ (cf. 
Definition 1). Then $\mathrm{Pol}(x, \pi) \neq \emptyset$ 
if and only if $S(\pi) = \mathbf{d}.$
\end{Thm}
 
{\em Proof}. The proof of this theorem can be found 
in \cite{He}, Theorem 7.1, (a). But we prove here that 
$\mathrm{Pol}(x, \pi) \neq \emptyset$ if $S(\pi) = \mathbf{d}$ 
because we will later use this argument. 
There is a basis $\{e(i,j)\}$ of $V$ indexed by the Young diagram 
$\mathbf{d}$ with the following properties (cf. \cite{S-S}, p.259, 
see also \cite{C-M}, 5.1.) 
\vspace{0.12cm}
 
(i) $\{e(i,j)\}$ is a Jordan basis of $x$, that is, $xe(i,j) = e(i-1,j)$ for 
$(i,j) \in \mathbf{d}$. 
\vspace{0.12cm}
 
(ii) $<e(i,j), e(p,q)> \neq 0$ if and only if $p = d_j -i + 1$ and 
$q = \beta(j)$, where $\beta$ is a permutation of 
$\{1, 2, ..., d^1\}$ which satiesfies:   
$\beta^2 = id$, $d_{\beta(j)} = d_j$, and 
$\beta(j) \not\equiv j$ (mod 2)
if $d_j \not\equiv \epsilon$ (mod 2). One can choose an arbitrary 
$\beta$ within these restrictions.   
\vspace{0.12cm}
 
For a sequence $(p_1, ..., p_s)$ with 
$\pi = \mathrm{ord}(p_1, ..., p_s)$ 
and $p_i = p_{s+1-i}, (1 \leq i \leq s)$, 
we shall construct an admissible 
flag $F$ of type $(p_1, ..., p_s)$ such that $xF_i \subset F_{i-1}$ 
for all $i$. We proceed by the induction on $s$. When 
$s = 1$, $\pi = [1^n]$ and $\pi = \mathbf{d}$. In this case, 
$x = 0$ and $F$ is a trivial flag $F_1 = V$. 
When $s > 1$, we shall construct an isotropic flag 
$0 \subset F_1 \subset F_{s-1} \subset V$. Put $p 
:= p_1(=p_s)$ and let $\rho := \mathrm{ord}(p_2, ..., p_{s-1}) 
\in \mathrm{Pai}(n-2p, q).$ Then we have 
 
$$\rho_j := \left\{ 
\begin{array}{rl}    
\pi_j -2 & \quad (j \leq p)\\ 
\pi_j  & \quad (j > p) 
\end{array}\right.$$ 

Let $$S': \mathrm{Pai}(n-2p, q) \to 
P_{\epsilon}(n-2p)$$ be the Spatenstein 
map and we put $\mu := S'(\rho).$ 
There are two cases (A) and (B). The first case 
(A) is when $i(\pi) = \{p\} \cup I(\rho)$ and 
$p \notin I(\rho)$. In this case, 
$p \not\equiv n$ (mod 2), $\pi_p \equiv \epsilon$ 
(mod 2) and $\pi_p = \pi_{p+1}+2$. Now we have 
$$ \mu_j = d_j -2, (j < p),$$  
$$ \mu_p = d_p -1, $$ 
$$ \mu_{p+1} = d_{p+1} -1, $$ 
$$ \mu_j = d_j, (j > p+1),$$ 
where $d_p = d_{p+1}$.  
The second case is exactly when (A) does not 
occur. In this case, $I(\pi) = I(\rho)$ and 
$$ \mu_j = d_j -2, (j \leq p),$$ 
$$ \mu_j = d_j, (j > p).$$   
Let us assume that the case (A) occurs. 
We choose the basis $e(i,j)$ of $V$ in such a way that 
the permutaion $\beta$ satisfies $\beta (p) = p+1$. 
There are two choices for $F_1$. The first one is to put 
$$ F_1 = \Sigma_{1 \leq j \leq p}\mathbf{C}e(1,j). $$  
The second one is to put 
$$ F_1 = \Sigma_{1 \leq j \leq p+1, j \neq p}
\mathbf{C}e(1,j). $$   
In any case, we put $F_{s-1} = F_1^{\perp}$. 
Then $x$ induces a nilpotent endomorphism of 
$F_{s-1}/F_1$ of type $\mu$. 
Next assume that the case (B) occurs. In this case, 
we put 
$$ F_1 =  \Sigma_{1 \leq j \leq p}\mathbf{C}e(1,j)$$ 
and $F_{s-1} = F_1^{\perp}$. Then $x$ induces a 
nilpotent endomorphism of $F_{s-1}/F_1$ of type $\mu$. 
By the induction on $s$, we have an admissible filtration 
$0 \subset F_1 \subset ... \subset F_{s-1} \subset V$ 
with desired properties.  
Let $P$ be the stabilizer group of $F$. Then 
it is clear that $x \in \mathfrak{n}(P)$. By an explicit calculation 
$\dim \mathcal{O}_x = 2 \dim G/P$. 

\begin{Thm}
Let $G$ and $\g$ be the same as Theorem 2.7. 
Let $x \in \g$ be a nilpotent 
element of type $\mathbf{d}$ 
and denote by $\mathcal{O}$ the orbit containing 
$x$. Assume that $P$ is a polarization of $x$ 
with Levi type $\pi$. Let 
$$\mu : T^*(G/P) \to \bar{\mathcal{O}}$$ 
be the Springer map. Then 
$$\mathrm{deg}(\mu) := \left\{ 
\begin{array}{rl}    
2^{\sharp I(\pi)-1}& \quad (q = \epsilon = 0, 
\pi^i \not\equiv 0\,(\mathrm{mod}\,2)\,\exists i)\\ 
2^{\sharp I(\pi)}  & \quad (q + \epsilon 
\geq 1\: \mathrm{or}\, q = \epsilon = 0,\, 
\pi^i \equiv 0\,(\mathrm{mod}\,2)\,\forall i) 
\end{array}\right.$$ 
Moreover, if $\mathrm{deg}(\mu) = 1$, 
then the Levi type of $P$ is 
unique. In other words, if two 
polarizations of $x$ respectively give Springer 
resolutions of $\bar{\mathcal{O}}$, 
then they have the same Levi type. 
\end{Thm}

{\em Proof}. The first part is \cite{He}, Theorem 
7.1, (d) (cf. \cite{He}, \S 1). 
The proof of the second part is rather 
technical, but for the completeness, 
we include it here. Let 
$$B(\mathbf{d}) = \{j \in \mathbf{N}; d_j > d_{j+1}, 
d_j \not\equiv \epsilon \,(\mathrm{mod}\,2)\}.$$ 
Note that $S(\pi) = \mathbf{d}$, where $S$ is the 
Spaltenstein map. When $\epsilon = 0$, 
$B(\mathbf{d}) = \emptyset$ if and only if 
$q = 0$ and $d^i \equiv 0$ (mod 2) for 
all $i$. 
Assume that $B(\mathbf{d}) = \emptyset$. 
Since $\mathrm{deg}(\mu) = 1$, 
by the first part of our theorem, 
$\sharp I(\pi) = 0$.  
Then $\pi = \mathbf{d}$.   
Assume that $B(\mathbf{d}) \neq 
\emptyset$. If $q \neq 0$ for our $\pi$ 
or $\epsilon = 1$, then $\sharp I(\pi) =  0$; 
hence $\pi = \mathbf{d}$. If $\epsilon 
= 0$ and $q = 0$ for $\pi$, then 
$\sharp I(\pi) = 1$. Since $\sharp I(\pi) 
= 1/2\,\sharp \{j; d_j \equiv 1\,(\mathrm{mod}\, 2)\}$ 
by \cite{He}, Lemma 6.3, (b). This implies 
that $\sharp\{j; d_j \equiv 1\,(\mathrm{mod}\,2)\} = 2$. 
Note that $\pi$ with $q = 0$ is uniquely 
determined by $\mathbf{d}$ because the 
Spaltenstein map is injective (\cite{He}, 
Prop. 6.5, (a)).  

Now let us prove the second part of 
our theorem. 
When $\epsilon = 1$, we should 
have $\pi = \mathbf{d}$ by the argument 
above. 
Next consider the case where $\epsilon = 0$. 
Assume that 
there exist two polarizations $P_1$ 
and $P_2$ giving Springer resolutions. 
Let $\pi_1$ and $\pi_2$ be their 
Levi types. Assume that 
$\pi_1 \in \mathrm{Pai}(n,0)$ 
and $\pi_2 \in \mathrm{Pai}(n,q_2)$ 
with $q_2 > 0$. 
By the argument above, we see that 
$\sharp\{j; d_j \equiv 1\,(\mathrm{mod}\,2)\} = 2$. 
On the other hand, since $q_2 > 0$, 
$\pi_2 = \mathbf{d}$. This shows  
that $q_2 = 2$; but, when $\epsilon = 0$, 
$q_2 \neq 2$ by Example 2.3, 
which is a contradiction.      
Hence, in this case, $\pi$ is 
also uniquely determined by $\mathbf{d}$. 
   
\section{Mukai flops}

{\bf Mukai flop of type A}.  
Let $x \in \mathfrak{sl}(n)$ be a 
nilpotent element of type $[2^k, 1^{n-2k}]$ and 
let $\mathcal{O}$ be the nilpotent orbit containing 
$x$. By Theorem 2.6, there are two polarizations 
$P$ and $P'$ of $x$, 
where $P$ has the flag type $(k, n-k)$ and $P'$ has 
the flag type $(n-k,k)$.  
The closure $\bar{\mathcal{O}}$ of $\mathcal{O}$ admits 
two Springer resolutions 
$$ T^*(SL(n)/P) \stackrel{\pi}\rightarrow \bar{\mathcal{O}} 
\stackrel{\pi'}\leftarrow T^*(SL(n)/P'). $$ 
Note that $SL(n)/P$ is isomorphic to the Grassmannian 
$G(k,n)$ and 
$SL(n)/P'$ is isomorphic to $G(n-k,n)$. 

\begin{Lem}
When $k < n/2$, $\pi$ and $\pi'$ are both small birational 
maps and the diagram becomes a flop. 
\end{Lem}

{\em Proof}. The closure $\bar{\mathcal{O}}$ consists of 
finite number of orbits 
$\{\mathcal{O}_{[2^i,1^{n-2i}]}\}_{0 \leq i \leq k}$. 
The main orbit $\mathcal{O}_{[2^k, 1^{n-2k}]}$ is an 
open set of $\bar{\mathcal{O}}$. A fiber of $\pi$ (resp. 
$\pi'$) over a point of $\mathcal{O}_{[2^i, 1^{n-2i}]}$ is 
isomorphic to the Grassmannian $G(k-i, n-2i)$ (resp. 
$G(n-i-k, n-2i)$). By a simple dimension count, if 
$k < n/2$, then $\pi$ 
and $\pi'$ are both small birational maps. 
Next let us prove that the diagram is a flop. 
Let $\tau \subset \mathcal{O}^{\oplus n}_{G(k,n)}$ (resp. 
$\tau' \subset \mathcal{O}^{\oplus n}_{G(n-k,n)}$) be the 
universal subbundle. Denote by $T$ (resp. $T'$) the pull-back 
of $\tau$ (resp. $\tau'$) by the projection 
$T^*G(k,n) \to G(k,n)$ (resp. $T^*G(n-k,n) \to G(n-k,n)$).  
We shall describe the strict transform of $\wedge^k 
T$ by the birational map $T^*G(k,n) --\to T^*G(n-k,n)$. 
Take a point $y \in \mathcal{O}_{[2^k, 1^{n-2k}]}$. 
Note that $T^*G(k,n)$ is  
naturally embedded in $G(k,n) \times \bar{\mathcal{O}}$. 
Then the fiber $\pi^{-1}(y)$ consists of one point 
$([\mathrm{Im}(y)], y) \in G(k,n) \times \bar{\mathcal{O}}.$ 
The fiber $T_{\pi^{-1}(y)}$ of the vector bundle $T$ over 
$\pi^{-1}(y)$ coincides with the vector space 
$\mathrm{Im}(y)$. Hence $(\wedge^k T)_{\pi^{-1}(y)}$ is 
isomorphic to $\wedge^k \mathrm{Im}(y)$. 
Now let $L$ be the strict transform of $\wedge^k T$ by 
$T^*G(k,n) --\to T^*G(n-k,n)$.  
First note that $(\pi')^{-1}(y)$ also consists of one 
point $([\mathrm{Ker}(y)], y) \in G(n-k,n) \times 
\bar{\mathcal{O}}$. Then, by definition, $L_{(\pi')^{-1}(y)} 
= \wedge^k \mathrm{Im}(y)$. Since 
$\wedge^k \mathrm{Im}(y) \cong (\wedge^{n-k}\mathrm{Ker}(y))^*$, 
we see that $L \cong (\wedge^{n-k} T')^{-1}$.        
Now $\wedge^k T$ is a negative line bundle. 
On the other hand, its strict transform $L$ becomes an 
ample line bundle. This implies that our diagram is 
a flop. 

\begin{Rque}
When $k = n/2$, $\pi$ and $\pi'$ are both 
divisorial birational contraction maps. 
Moreover, two resolutions are isomorphic. 
\end{Rque}

\begin{Def}
The diagram 
$$ T^*(SL(n)/P) \stackrel{\pi}\rightarrow \bar{\mathcal{O}} 
\stackrel{\pi'}\leftarrow T^*(SL(n)/P') $$ 
is called a (stratified) Mukai flop of type A when $k < n/2$. 
\end{Def}
\vspace{0.12cm}
{\bf Mukai flop of type D}. 
Assume that $k$ is an odd integer 
with $k \geq 3$. Let $V$ be a 
$\mathbf{C}$-vector space of $\dim 2k$ with a 
non-degenerate symmetric form $<,>$.  
Let $x \in \mathfrak{so}(V)$ be a nilpotent 
element of type $[2^{k-1}, 1^2]$ and let 
$\mathcal{O}$ be the nilpotent orbit 
containing $x$. 
Let $S: \mathrm{Pai}(2k,0) \to P_{\epsilon}(2k)$ 
be the Spaltenstein map, where $\epsilon = 0$ 
in our case. Then, for $\pi := (2^k) \in 
\mathrm{Pai}(2k,0)$, $S(\pi) = [2^{k-1}, 1^2]$. 
Let us recall the construction of the stabilized 
flags by the polarizations of $x$ in the proof of 
Theorem 2.7. Since $I(\pi) = \{k\}$, the case 
(A) occurs (cf. the proof of Theorem 2.7.); hence 
there are two choices of the flags. We denote by 
$P^+$ the stabilizer subgroup of $SO(V)$ of one 
flag, and denote by $P^-$ the stabilizer subgroup 
of another one. Let $G_{iso}(k,V)$ be the orthogonal 
Grassmannian which parametrizes $k$ dimensional 
isotropic subspaces of $V$. $G_{iso}(k, V)$ has 
two connected components ${G^+}_{iso}(k,V)$ and 
${G^-}_{iso}(k,V)$. Note that 
$SO(V)/P^+ \cong {G^+}_{iso}(k,V)$ and 
$SO(V)/P^- \cong {G^-}_{iso}(k,V)$.    
The closure $\bar{\mathcal{O}}$ of $\mathcal{O}$ 
admits two Springer resolutions 
$$ T^*(SO(V)/P^+) \stackrel{\pi^+}\rightarrow 
\bar{\mathcal{O}} \stackrel{\pi^-}\leftarrow 
T^*(SO(V)/P^-). $$ 

\begin{Lem}
$\pi^+$ and $\pi^-$ are both small birational maps and 
the diagram becomes a flop. 
\end{Lem}

{\em Proof}. The closure $\bar{\mathcal{O}}$ consists 
of the orbits 
$\{\mathcal{O}_{[2^{k-2i-1}, 1^{4i+2}]}\}_{1 
\leq i \leq 1/2(k-1)}$. The main orbit is an open 
set of $\bar{\mathcal{O}}$. A fiber of 
$\pi^+$ (resp. $\pi^-$) over a point of 
$\mathcal{O}_{[2^{k-2i-1}, 1^{4i+2}]}$ is 
isomorphic to ${G^+}_{iso}(2i+1,4i+2)$ 
(resp. ${G^-}_{iso}(2i+1,4i+2)$). By dimension 
counts of each orbit and of each fiber, we see that 
$\pi^+$ and $\pi^-$ are both small birational maps.  
Next let us prove that the diagram is a flop. 
Let $\tau^+ \subset \mathcal{O}^{\oplus 2k}_{{G^+}_{iso}(k,V)}$ 
(resp. 
$\tau^- \subset \mathcal{O}^{\oplus 2k}_{{G^-}_{iso}(k,V)}$) be 
the universal subbundle. Denote by $T^+$ (resp. 
$T^-$)  the pull-back of $\tau^+$ (resp. 
$\tau^-$) by the projection $T^*({G^+}_{iso}(k,V)) 
\to {G^+}_{iso}(k,V)$ (resp. 
$T^*({G^-}_{iso}(k,V)) \to {G^-}_{iso}(k,V)$). 
We shall describe the strict transform of 
$\wedge^k T^-$ by the birational map 
$T^*({G^-}_{iso}(k,V)) --\to 
T^*({G^+}_{iso}(k,V))$. Take a point $y 
\in \mathcal{O}_{[2^{k-1}, 1^2]}$. 
Let $g \in SO(V)$ be an element such that 
$gxg^{-1} = y$. Note that 
$T^*({G^+}_{iso}(k,V))$ (resp. $T^*({G^-}_{iso}(k,V))$) 
is naturally embedded 
in ${G^+}_{iso}(k,V) \times \bar{\mathcal{O}}$ 
(resp. ${G^-}_{iso}(k,V) \times \bar{\mathcal{O}}$). 
Then the fiber $(\pi^+)^{-1}(y)$ (resp. 
$(\pi^-)^{-1}(y)$) 
consists 
of one point $([F^+_y], y) \in  
{G^+}_{iso}(k,V) \times \bar{\mathcal{O}}$ 
(resp. $([F^-_y], y) \in  
{G^-}_{iso}(k,V) \times \bar{\mathcal{O}}$)   
where $F^+_y \subset V$ (resp. 
$F^-_y \subset V$) 
is the flag stabilized by $gP^+g^{-1}$ (resp.  
$gP^-g^{-1}$).  Note that 
$gP^+g^{-1}$ and $gP^-g^{-1}$ are both polarizations of $y$.  
Let us recall the construction of flags in 
the proof of Theorem 2.7. 
For $y$ we choose a Jordan basis $\{e(i,j)\}$ 
of $V$ as in the proof of Theorem 2.7. 
Since $\mathbf{d} = [2^{k-1}, 1^2]$, 
$\beta$ is a permutation of 
$\{1, 2, ..., k, k+1\}$. But it preserves 
the subsets $\{1, 2, ..., k-1\}$ and $\{k,k+1\}$ 
respectively. We assume that $\beta (k) 
= k+1$ and $\beta (k+1) = k$. 
In our situation, the case (A) occurs. There 
are two choices of the flags: 
$$\Sigma_{1 \leq j 
\leq k-1}\mathbf{C}e(1,j) + \mathbf{C}e(1,k)$$ 
and  
$$\Sigma_{1 \leq j 
\leq k-1}\mathbf{C}e(1,j) + \mathbf{C}e(1,k+1).$$ 
Note that one of these is stabilized by 
$gP^+g^{-1}$ and another one is stabilized 
by $gP^-g^{^-1}$. We may assume that   
$$ F^+_y = \Sigma_{1 \leq j 
\leq k-1}\mathbf{C}e(1,j) + \mathbf{C}e(1,k), $$ 
and  
$$ F^-_y = \Sigma_{1 \leq j 
\leq k-1}\mathbf{C}e(1,j) + \mathbf{C}e(1,k+1).$$ 
Since $\mathrm{Ker}(y) = \Sigma_{1 \leq j \leq k+1}e(1,j)$ 
and $\mathrm{Im}(y) = \Sigma_{1 \leq j \leq k-1}e(1,j)$, 
we have two exact sequences 
$$ 0 \to \mathrm{Ker}(y)/F^+_y  \to V/F^+_y  
\to \mathrm{Im}(y) \to 0, $$
and 
$$ 0 \to \mathrm{Im}(y) \to F^-_y \to F^-_y/\mathrm{Im}(y) 
\to 0.$$ 
Since $F^-_y/\mathrm{Im}(y) \cong \mathrm{Ker}(y)/F^+_y$, 
we conclude that 
$$ \wedge^k F^-_y \cong \wedge^k (V/F^+_y).$$ 
Let $L$ be the strict transform of $\wedge^k T^-$ by the 
birational map $T^*({G^-}_{iso}(k,V)) --\to 
T^*({G^+}_{iso}(k,V))$. The fiber $T^-_{(\pi^-)^{-1}(y)}$ 
of the vector bundle $T^-$ is isomorphic to the vector 
space $\wedge^k F^-_y$. Hence, by the definition of $L$, 
$L_{(\pi^+)^{-1}(y)} = \wedge^k F^-_y$. By the observation 
above, we see that 
$L_{(\pi^+)^{-1}(y)} = \wedge^k (V/F^+_y)$. This shows 
that $L \cong (\wedge^k T^+)^{-1}$. 
Now $\wedge^k T^-$ is a negative line bundle. On the other 
hand, its strict transform $L$ is an ample line bundle. 
This implies that our diagram is a flop. 

\begin{Def}
The diagram        
$$ T^*(SO(V)/P^+) \stackrel{\pi^+}\rightarrow 
\bar{\mathcal{O}} \stackrel{\pi^-}\leftarrow 
T^*(SO(V)/P^-)$$ 
is called a (stratified) Mukai flop of type 
$D$ or called an orthogonal Mukai flop. 
\end{Def}

\begin{Rque}
When $k$ is an even integer with $k \geq 2$, 
there are two nilpotent orbits $\mathcal{O}^+$ 
and $\mathcal{O}^-$ with Jordan type $[2^k]$. 
They have Springer resolutions 
$$ T^*({G^+}_{iso}(k,2k)) \to \bar{\mathcal{O}}^+, $$ 
and $$ T^*({G^-}_{iso}(k,2k)) \to 
\bar{\mathcal{O}}^-.$$ 
These resolutions are both divisorial birational 
contraction maps.  
When $k = 1$, three varieties 
$T^*({G^+}_{iso}(1,2))$, $T^*({G^-}_{iso}(1,2))$ 
and $\bar{\mathcal{O}}$ are all isomorphic.  
\end{Rque}

\begin{Def}
Let $X \stackrel{f}\rightarrow Y 
\stackrel{f'}\leftarrow X'$ be two resolutions 
of a variety $Y$. Then this diagram is called 
a locally trivial family of Mukai flops 
of type A (resp. of type D) if there is a 
partial open covering $\{U_{\lambda}\}$ of 
$Y$ which contains the singular locus of 
$Y$ such that each diagram 
$$f^{-1}(U_{\lambda}) \rightarrow 
U_{\lambda} \leftarrow (f')^{-1}(U_{\lambda})$$ 
is isomorphic to the product of a Mukai flop 
of type A (resp. of type D) with a suitable 
disc $\Delta^m$.   
\end{Def}

\section{Symplectic resolutions of Nilpotent 
orbits}

In this section, we shall prove that any two 
symplectic resolutions of a nilpotent orbit 
closure of a classical simple Lie algebras 
are connected by a finite sequence of diagrams  
which are locally trivial families of 
Mukai flops.
 
\begin{Lem} 
Let $x \in \mathfrak{sl}(n)$ be a 
nilpotent element of type 
$\mathbf{d} := [d_1, ..., d_k]$. 
Let $(p_1, ..., p_s)$ be a sequence 
of positive integers such that 
$\mathrm{ord}(p_1, ..., p_s) 
= \mathbf{d}$. Fix a flag 
$F := \{F_i\}$ of $V := \mathbf{C}^n$ 
of type $(p_1, ..., p_s)$ such that 
$xF_i \subset F_{i-1}$ for all $i$. 
Assume that $p_j \neq p_{j-1}$ for 
an index $j$. Then we obtain a new 
flag $F'$ of type $(p_1, ..., 
p_j, p_{j-1}, ..., p_s)$ from $F$ such that 
$xF'_i \subset F'_{i-1}$ for all $i$ by 
the following operation. 
\vspace{0.12cm}

(The case where $p_{j-1} < p_j$): 
$x$ induces an endomorphism $\bar{x} 
\in \mathrm{End}(F_j/F_{j-2})$. For 
the projection $\phi: F_j \to 
F_j/F_{j-2}$, we put $F'_{j-1} := 
\phi^{-1}(\mathrm{Ker}(\bar{x}))$. We 
then put 
$$F'_i := \left\{ 
\begin{array}{rl}    
F_i & \quad (i \neq j-1)\\ 
F'_{j-1} & \quad (i = j-1) 
\end{array}\right.$$

(The case where $p_{j-1} > p_j$): 
$x$ induces an endomorphism $\bar{x} 
\in \mathrm{End}(F_j/F_{j-2})$. For 
the projection $\phi: F_j \to 
F_j/F_{j-2}$, we put $F'_{j-1} := 
\phi^{-1}(\mathrm{Im}(\bar{x}))$. We 
then put 
$$F'_i := \left\{ 
\begin{array}{rl}    
F_i & \quad (i \neq j-1)\\ 
F'_{j-1} & \quad (i = j-1) 
\end{array}\right.$$
\end{Lem}
\vspace{0.12cm}

\begin{Lem}
Let $V$ be a  
$\mathbf{C}$-vector space of $\dim n$ 
with a non-degenerate 
bilinear form such that $<v,w> = (-1)^{\epsilon}
<w,v>$ for all $v,w \in V$. Let $\g$ be the Lie 
algebra $\mathfrak{so}(V)$ or $\mathfrak{sp}(V)$ 
according as $\epsilon = 0$ or $\epsilon = 1$. 
Let $x\in \g$ be a nilpotent element of type 
$\mathbf{d}$. Suppose that for $\pi \in 
\mathrm{Pai}(n,q)$, $\mathbf{d} = S(\pi)$ where 
$S$ is the Spaltenstein map. Let $(p_1, ..., p_k, 
q, p_k, ..., p_1)$ be a sequence of integers such 
that $\pi = \mathrm{ord}(p_1, ..., p_k,q, p_k, ..., 
p_1)$. Fix an admissible flag $F$ of type 
$(p_1, ..., p_k,q,p_k, ..., p_1)$ such that 
$xF_i \subset F_{i-1}$ for all $i$. 
\vspace{0.12cm}

(i) Assume that $p_{j-1} \neq p_j$ for an index 
$1 \leq j \leq k$.   
Then we obtain a new 
flag $F'$ of type $(p_1, ..., 
p_j, p_{j-1}, ..., 
p_k,q,p_k, ..., p_{j-1},p_j, ..., p_1)$ from 
$F$ such that 
$xF'_i \subset F'_{i-1}$ for all $i$ by 
the following operation. 
\vspace{0.12cm}

(The case where $p_{j-1} < p_j$): 
$x$ induces an endomorphism $\bar{x} 
\in \mathrm{End}(F_j/F_{j-2})$. For 
the projection $\phi: F_j \to 
F_j/F_{j-2}$, we put $F'_{j-1} := 
\phi^{-1}(\mathrm{Ker}(\bar{x}))$. We 
then put 
$$F'_i := \left\{ 
\begin{array}{rl}    
F_i & \quad (i \neq j-1,2k+2-j)\\ 
F'_{j-1} & \quad (i = j-1)\\
(F'_{j-1})^{\perp} & \quad (i = 2k+2-j) 
\end{array}\right.$$

(The case where $p_{j-1} > p_j$): 
$x$ induces an endomorphism $\bar{x} 
\in \mathrm{End}(F_j/F_{j-2})$. For 
the projection $\phi: F_j \to 
F_j/F_{j-2}$, we put $F'_{j-1} := 
\phi^{-1}(\mathrm{Im}(\bar{x}))$. We 
then put 
$$F'_i := \left\{ 
\begin{array}{rl}    
F_i & \quad (i \neq j-1,2k+2-j)\\ 
F'_{j-1} & \quad (i = j-1)\\
(F'_{j-1})^{\perp} & \quad (i = 2k+2-j) 
\end{array}\right.$$

(ii) Assume that $q = 0$ and $p_k$ is 
odd. Then there is an admissible flag $F'$ 
of $V$ of type $(p_1, ..., p_k,p_k, ...,p_1)$ 
such that 
\vspace{0.12cm}

$xF'_i \subset F'_{i-1}$ for all $i$,  
\vspace{0.12cm}
  
$F'_i =     
F_i$ for $i \neq k$ 
and 
\vspace{0.12cm}

$F'_k \neq F_k$.  
\end{Lem}

{\em Proof}. (i): When $p_{j-1} 
< p_j$, $\mathrm{rank}(\bar{x}) = p_{j-1}$ for 
$\bar{x} \in \mathrm{End}(F_j/F_{j-2})$. 
In fact, since $xF_j \subset F_{j-1}$, 
$\mathrm{rank}(\bar{x}) \leq p_{j-1}$. 
Assume that $\mathrm{rank}(\bar{x}) 
< p_{j-1}$. Then we can construct a new 
flag from $F$ by replacing $F_{j-1}$ with 
 a subspace $F'_{j-1}$ containing $F_{j-2}$ 
such that 
$$ \mathrm{Im}(\bar{x}) \subset 
F'_{j-1}/F_{j-2} \subset \mathrm{Ker}(\bar{x})$$ 
and $\dim F'_{j-1}/F_{j-2} = p_{j-1}$. The new 
flag satisfies $xF'_i \subset F'_{i-1}$ for all $i$ 
and it has the same flag type as $F$. Since there 
are infinitely many choices of $F'_{j-1}$, we have 
infinitely many such $F'$. This contradicts the fact 
that $x$ has only finite polarizations. 
Hence, $\mathrm{rank}(\bar{x}) = p_{j-1}$. Then the 
flag $F'$ in our Lemma satisfies the desired 
properties. When $p_{j-1} > p_j$, we see that 
$\dim \mathrm{Ker}(\bar{x}) = p_{j-1}$ by a similar 
way. Then the latter argument is the same as 
when $p_{j-1} < p_j$. 
\vspace{0.12cm}

(ii): According to the proof of Theorem 2.7. we 
construct a flag $F$ such that $xF_i \subset F_{i-1}$. 
Since $q = 0$ and $p_k$ is odd, we have the case (A) 
in the last step. As a consequence, we have two 
choices of the flags. One of them is $F$ and another 
one is $F'$. \hspace{1.2cm} Q.E.D. 
\vspace{0.2cm}

Let $F$ be the flag in Lemma 4.1, Lemma 4.2, (i) 
or Lemma 4.2, (ii). In each lemma, we have 
constructed another flag $F'$. Let $G$ be the 
complex Lie group $SL(V)$, $Sp(V)$ or $SO(V)$ 
according as $V$ is simply a $\mathbf{C}$-vector 
space with no bilinear forms, with a non-degenerate 
skew-symmetric form or with a non-degenerate 
symmetric form. 
Let $P \subset G$ (resp. $P' \subset G$) 
be the stabilizer group of 
the flag $F$ (resp. $F'$). Then $P$ and 
$P'$ are both polarizations of $x \in \g$. 
Let $\mathcal{O} \subset \g$ be the nilpotent 
orbit containing $x$. Let us consider two 
Springer resolutions 
$$ T^*(G/P) \stackrel{\mu}\rightarrow \bar{\mathcal{O}} 
\stackrel{\mu'}\leftarrow T^*(G/P'). $$ 
Note that $T^*(G/P)$ (resp. $T^*(G/P')$) is 
embedded in $G/P \times \bar{\mathcal{O}}$ 
(resp. $G/P' \times \bar{\mathcal{O}}$). 
An element $y \in \mathcal{O}$ can be written 
as $y = gxg^{-1}$ for some $g \in G$. 
Then we have 
$$  \mu^{-1}(y) = ([gF], y) \in G/P 
\times \bar{\mathcal{O}}, $$ 
and 
$$  (\mu')^{-1}(y) = ([gF'], y) \in 
G/P' \times \bar{\mathcal{O}}. $$ 
We define the flag $\bar{F}$ in the following 
manner. If $F$ is the flag in Lemma 4.1, then 
$\bar{F}$ is the flag obtained from $F$ by 
deleting the subspace $F_{j-1}$. 
If $F$ is the flag in Lemma 4.2, (i), then 
$\bar{F}$ is the flag obtained from $F$ by 
deleting subspaces $F_{j-1}$ and 
$F_{2k+2-j}$. Finally, if $F$ is the 
flag in Lemma 4.2, (ii), then $\bar{F}$ is 
the flag obtained from $F$ by deleting 
$F_k$. Note that $\bar{F}$ is also obtained 
from $F'$ by the same manner. 
Let $\bar{P} \subset G$ be the stabilizer 
group of the flag $\bar{F}$. 
We then have two projections 
$$ G/P \stackrel{p}\rightarrow G/\bar{P} 
\stackrel{p'}\leftarrow G/P'. $$ 
By two projections  
$$ G/P \times \bar{\mathcal{O}} 
\stackrel{p \times id}\rightarrow 
G/\bar{P} \times \bar{\mathcal{O}} 
\stackrel{p' \times id}\leftarrow 
G/P' \times \bar{\mathcal{O}}, $$ 
$T^*(G/P)$ and $T^*(G/P')$ have 
the same image $X$ in  
$G/\bar{P} \times \bar{\mathcal{O}}$. 
Since $p$ and $p'$ are both proper maps, 
$X$ is a closed subvariety of 
$G/\bar{P} \times \bar{\mathcal{O}}$. 
The following diagram has been obtained 
as a consequence:
$$ T^*(G/P) \rightarrow X \leftarrow 
T^*(G/P'). $$ 

\begin{Lem}
When $F$ is the flag in Lemma 4.1 or 
in Lemma 4.2, (i), the diagram   
$$ T^*(G/P) \stackrel{f}\rightarrow X 
\stackrel{f'}\leftarrow 
T^*(G/P') $$ 
is locally a trivial family of  
Mukai flops of type A. 
When $F$ is the flag in Lemma 4.2, (ii), 
the diagram is locally a trivial 
family of Mukai flops of type D. 
\end{Lem}

{\em Proof}. We prove the assetion when 
$\g = \mathfrak{so}(n)$ or $\g = 
\mathfrak{sp}(n)$. The case when $\g 
= \mathfrak{sl}(n)$ is easier; so we omit 
the proof.  

Consider the situation in Lemma 4.2, (i). 
A point of $G/\bar{P}$ corresponds to an 
isotropic flag $\bar{F}$ of $V$ of type  
$(p_1, ..., p_{j-1}+p_j, ..., p_k, q, 
p_k, ..., p_{j-1}+p_j, ..., p_1)$. 
Let 
$$0 \subset \bar{\mathcal{F}}_1 \subset 
... \subset \bar{\mathcal{F}}_{2k-1} = 
(\mathcal{O}_{G/\bar{P}})^n$$ 
be the universal subundles on $G/\bar{P}$. 
Let $$W \subset 
\underline{\mathrm{End}}(\bar{\mathcal{F}}_{j-1}/
\bar{\mathcal{F}}_{j-2})$$ 
be the subvariety consisting of the points 
$([\bar{F}], \bar{x})$ where $\bar{x} \in 
\mathrm{End}(\bar{F}_{j-1}/\bar{F}_{j-2})$, 
$\bar{x}^2 = 0$ 
and $\mathrm{rank}(\bar{x}) \leq \mathrm{min}
(p_j, p_{j-1})$. If we put $m := p_{j-1} 
+ p_j$ and $r:= \mathrm{min}(p_j, p_{j-1})$, then 
$$W \to G/\bar{P}$$ is an 
$\bar{\mathcal O}_{[2^r, 1^{m-2r}]}$ 
bundle over $G/\bar{P}$. 
Let us recall the definition of $X$.    
$$X \subset G/\bar{P} \times \bar{\mathcal O}$$ 
consists of the points $([\bar{F}], x)$ such that 
$x\bar{F}_i \subset \bar{F}_{i-1}$ for all $i \neq 
j-1, 2k-j$ and $x\bar{F}_i \subset \bar{F}_i$ for 
$i = j-1, 2k-j$. Moreover, 
the induced endomorphism $\bar{x} \in 
\mathrm{End}(\bar{F}_{j-1}/\bar{F}_{j-2})$ satisfies 
$\bar{x}^2 
= 0$ and $\mathrm{rank}(\bar{x}) \leq 
\mathrm{min}(p_{j-1}, p_j)$. 
Let  
$$ \phi : X \to W $$ 
be the projection defined by 
$\phi([\bar{F}],x) = ([\bar{F}], \bar{x})$, 
where $\bar{x} \in 
\mathrm{End}(\bar{F}_{j-1}/\bar{F}_{j-2})$ 
is the induced endomorphism by $x$. 
It can be checked that $\phi$ is an affine 
bundle. Since $W$ is an  
$\bar{\mathcal O}_{[2^r, 1^{m-2r}]}$ bundle 
over $G/\bar{P}$, there exists a family of 
Mukai flops of type A: 
$$ Y \rightarrow W \leftarrow Y' $$ 
parametrized by $G/\bar{P}$. The diagram 
$$ T^*(G/P) \rightarrow X \leftarrow 
T^*(G/P') $$ 
coincides with the pull back of the 
previous diagram by $\phi: X \to W$. 
Since $\phi$ is an affine bundle, this diagram 
is locally a trivial family of Mukai flops 
of type A. 

Next consider the situation in Lemma 4.2, (ii). 
A point of $G/\bar{P}$ corresponds to an 
isotropic flag $\bar{F}$ of $V$ of type  
$(p_1, ..., 2p_k, ..., p_1)$. 
Let 
$$0 \subset \bar{\mathcal{F}}_1 \subset 
... \subset \bar{\mathcal{F}}_{2k-1} = 
(\mathcal{O}_{G/\bar{P}})^n$$ 
be the universal subundles on $G/\bar{P}$. 
Let $$W \subset \underline{\mathrm{End}}
(\bar{\mathcal{F}}_k/\bar{\mathcal{F}}_{k-1})$$ 
be the subvariety consisting of the points 
$([\bar{F}], \bar{x})$ where 
$$\bar{x} 
\in \bar{\mathcal{O}}_{[2^{p_k-1}, 1^2]} 
\subset \mathfrak{so}(\bar{F}_k/\bar{F}_{k-1}).$$  
$W \to G/\bar{P}$ is an 
$\bar{\mathcal O}_{[2^{p_k-1}, 1^2]}$ 
bundle over $G/\bar{P}$. 
Let us recall the definition of $X$.    
$$X \subset G/\bar{P} \times \bar{\mathcal O}$$ 
consists of the points $([\bar{F}], x)$ such that 
$x\bar{F}_i \subset \bar{F}_{i-1}$ for all $i \neq k$ 
and $x\bar{F}_k \subset \bar{F}_k$. Moreover, 
the induced endomorphism $\bar{x} \in 
\mathfrak{so}(\bar{F}_k/\bar{F}_{k-1})$ 
is contained in 
$\bar{\mathcal{O}}_{[2^{p_k-1}, 1^2]}$.
Let  
$$ \phi : X \to W $$ 
be the projection defined by 
$\phi([\bar{F}],x) = ([\bar{F}], \bar{x})$, 
where $\bar{x} \in 
\mathfrak{so}(\bar{F}_k/\bar{F}_{k-1})$ 
is the induced endomorphism by $x$. 
It can be checked that $\phi$ is an affine 
bundle. Since $W$ is an  
$\bar{\mathcal O}_{[2^{p_k-1}, 1^2]}$ bundle 
over $G/\bar{P}$, there exists a family of 
Mukai flops of type D: 
$$ Y \rightarrow W \leftarrow Y' $$ 
parametrized by $G/\bar{P}$. The diagram 
$$ T^*(G/P) \rightarrow X \leftarrow 
T^*(G/P') $$ 
coincides with the pull back of the 
previous diagram by $\phi: X \to W$. 
Since $\phi$ is an affine bundle, this diagram 
is locally a trivial family of Mukai flops 
of type D.

\begin{Thm}
Let $\mathcal{O} \subset \g$ be an  
orbit of a classical complex simple 
Lie algebra.  
Let $Y$ and $Y'$  
be any two Springer resolutions of 
the closure $\bar{\mathcal{O}}$ 
of the nilpotent orbit. 
Then the birational map $Y --\to Y'$ 
can be decomposed into finite number of 
diagrams $Y_i \rightarrow X_i \leftarrow Y_{i+1}$ 
with $Y_1 = Y$ and $Y_m = Y'$ in such a way 
that each diagram is locally a trivial family 
of Mukai flops of type $A$ or of type $D$. 
\end{Thm}

\begin{Rque}
By a theorem of Fu [Fu 1, Thm 3.3], 
any projective symplectic resolution of 
$\bar{\mathcal{O}}$ is 
obtained as a Springer resolution. 
\end{Rque}

{\em Proof}. (The case $\g = \mathfrak{sl}(V)$): 
We put $G = SL(V)$. Let $Y = T^*(G/P)$ 
and $Y' = T^*(G/P')$, where 
$P$ and $P'$ are polarizations of an 
element $x \in \mathcal{O}$. By Theorem 2.6, we may 
assume that $P$ has 
the flag type $(p_1, ..., p_s)$ 
and $P'$ has the flag type $(p_{\sigma(1)}, ..., 
p_{\sigma(s)})$ where $\sigma$ is a permutation 
of $\{1, 2, ..., s\}$. Let $F$ be the flag of $V$ 
stabilized by $P$. 
Applying the operations in Lemma 4.1 successively, 
one can reach a flag $F'$ of type $(p_{\sigma (1)}, ..., 
p_{\sigma (s)})$. Note that $P'$ is the stabilizer 
group of $F'$. By Lemma 4.3, each step corresponds 
to a diagram which is locally a trivial family 
of a Mukai flop of type A. 
\vspace{0.12cm}

(The case $\g = \mathfrak{sp}(V)$ or 
$\mathfrak{so}(V)$): Let $G$ be $SO(V)$ or $Sp(V)$. 
By Theorem 2.8, 
the parabolic subgroups giving Springer resolutions 
of $\bar{\mathcal{O}}$ all have the same Levi type 
$\pi \in \mathrm{Pai}(n,q)$. 
\vspace{0.12cm}

Assume that $q \neq 0$ or $\epsilon = 1$. 
We write 
$\pi = \mathrm{ord}(p_1, ..., p_k,q,p_k, ..., p_1)$ 
with a non-decreasing sequence $p_1 \leq p_2 \leq 
... \leq p_k$. 
Fix an element $x \in \mathcal{O}$. 
Let $Y = T^*(G/P)$ with a polarization  
$P$ of $x$.  The flag type of $P$ is 
$(p_{\sigma (1)}, ..., 
p_{\sigma (k)},q, p_{\sigma (k)}, ..., p_{\sigma (1)})$ 
with a permutation $\sigma$ of $\{1,2, ..., k\}$. 
Since $q \neq 0$ or $\epsilon = 1$, the conjugacy class 
of a parabolic subgroups is uniquely determined 
by its flag type. Let $P_0$ be a polarization of $x$  
with the flag type $(p_1, ..., p_k,q,p_k, ..., p_1)$. 
Let $F$ be an admissible flag of $V$ stabilized by 
$P$. Applying the operations in Lemma 4.2, (i) 
successively, one can reach a flag $F'$ of type 
$(p_1, ..., p_k,q,p_k, ..., p_1)$.  
Then $P_0$ is the stabilizer group of $F'$. 
By Lemma 4.3, each step corresponds 
to a diagram which is locally a trivial 
family of a Mukai flop 
of type A. Therefore, the birational map 
$T^*(G/P) --\to T^*(G/P_0)$ is decomposed into  
desired diagrams. 
By the same argument, the birational map 
$T^*(G/P') --\to T^*(G/P_0)$ is also decomposed 
into desired diagrams. This shows that   
$Y$ and $Y'$ are connected by diagrams which are 
locally trivial families of Mukai flops of 
type A.
\vspace{0.12cm}

Assume that $q = 0$ and $\epsilon = 0$.      
Since we are only concerned with parabolic 
subgroups whose Springer maps are birational, 
the members of ${}^t\pi$ (cf. Notation and Convention) 
are all even numbers or 
even if some odd numbers appear in ${}^t\pi$, they  
are all the same number (possibly with some 
multiplicity). 
\vspace{0.12cm}

(a) {\em The case where all members of 
${}^t\pi$ are even}:  

We write 
$\pi = \mathrm{ord}(p_1, ..., p_k,p_k, ..., p_1)$ 
with a non-decreasing sequence $p_1 \leq p_2 \leq 
... \leq p_k$.  
Our nilpotent orbit $\mathcal{O}$ is one of the 
nilpotent orbits with a {\em very even} partition. 
Fix an element $x \in \mathcal{O}$. Then a 
polarization of $x$ is uniquely determined by 
its flag type. Let $P_0$ be a polarization of 
$x$ with the flag type $(p_1, ..., p_k,p_k, ..., p_1)$.  
The symplectic resolution $Y$ can be written as 
$T^*(G/P)$ with a polarization $P$ of $x$. 
The flag type of $P$ is 
$(p_{\sigma (1)}, ..., 
p_{\sigma (k)}, p_{\sigma (k)}, ..., p_{\sigma (1)})$ 
with some permutation $\sigma$ of 
$\{1, 2, ..., k\}$. 
Let $F$ be an admissible flag of $V$ stabilized 
by $P$.
Applying the operations in Lemma 4.2, (i) 
successively, one can reach a flag $F'$ of type 
$(p_1, ..., p_k,p_k, ..., p_1)$.  
Then $P_0$ is the stabilizer group of $F'$. 
By Lemma 4.3, each step corresponds 
to a diagram which is locally a trivial 
family of a Mukai flop 
of type A. Therefore, the birational map 
$T^*(G/P) --\to T^*(G/P_0)$ is decomposed into  
desired diagrams. 
On the other hand, another symplectic resolution 
$Y'$ is written as $T^*(G/P')$ with a polarization 
$P'$ of $x$. 
By the same argument, the birational map 
$T^*(G/P') --\to T^*(G/P_0)$ is also decomposed 
into desired diagrams. This shows that   
$Y$ and $Y'$ are connected by diagrams which are 
locally trivial families of Mukai flops of 
type A.  
\vspace{0.12cm}

(b) {\em The case where odd numbers appear 
in ${}^t\pi$}: 

One can write $\pi = 
\mathrm{ord}(p_1, ..., p_k,p_k, ..., p_1)$ 
in such a way that $p_1 \leq ... \leq p_l$ 
are even numbers and that $p_{l+1} = ... = 
p_k$ are odd. Fix $x \in \mathcal{O}$, and 
write $Y = T^*(G/P)$ and $Y' = T^*(G/P')$ 
with polarizations $P$ and $P'$ of $x$. 
Let $(p_{\sigma (1)}, ..., p_{\sigma (k)}, 
p_{\sigma (k)}, ..., p_{\sigma (1)})$ 
and $(p_{\tau (1)}, ..., p_{\tau (k)}, 
p_{\tau (k)}, ..., p_{\tau (1)})$ be the 
flag types of $P$ and $P'$ respectively. 
Let $F$ (resp. $F'$) be an admissible 
flag stabilized 
by $P$ (resp. $P'$). 
Applying the operations in Lemma 4.2, (i) 
to $F$ (resp. $F'$) 
successively, one can reach a flag $F^0$ 
(resp. $(F')^0$) of type 
$(p_1, ..., p_k,p_k, ..., p_1)$.
Let $P_0$ (resp. $P'_0$) be the stabilizer 
group of $F^0$ (resp. $(F')^0$). 
Then the birational maps  
$T^*(G/P) --\to T^*(G/P_0)$ and 
$T^*(G/P') --\to T^*(G/(P')_0)$ are 
decomposed into diagrams which are locally 
trivial families of Mukai flops of type A. 
If $p_{l+1} = ... = p_k = 1$, then 
$P_0$ and $P'_0$ are conjugate to each other. 
In this case, two Springer resolutions $T^*(G/P_0) 
\to \bar{\mathcal{O}}$ and $T^*(G/(P')_0) 
\to \bar{\mathcal{O}}$ are isomorphic. 
Then the birational map 
$T^*(G/P) --\to T^*(G/P')$ is decomposed 
into diagrams which are locally trivial 
families of Mukai flops of type A. 
If $p_{l+1} = ... = p_k > 1$, then 
$P_0$ and $P'_0$ may not be conjugate. 
When $P_0$ and $P'_0$ are conjugate, 
the birational map 
$T^*(G/P) --\to T^*(G/P')$ is decomposed 
into diagrams which are locally trivial 
families of Mukai flops of type A. 
When $P_0$ and $P'_0$ are not conjugate, 
the relationship between two flags $F^0$ 
and $(F')^0$ are described in Lemma 4.2, (ii). 
In this case, $T^*(G/P_0) --\to T^*(G/P'_0)$ 
is locally a trivial family of Mukai flops 
of type D. Then the birational map 
$T^*(G/P) --\to T^*(G/P')$ is decomposed 
into diagrams which are locally trivial families 
of Mukai flops of type A and of type D.   

\begin{Exam}
Let $\mathcal{O} \subset \mathfrak{sl}(6)$ 
be the nilpotent orbit of Jordan type 
$[3,2,1]$. Take an element $x \in 
\mathcal{O}$. Then $x$ has six polarizations 
$P_{\sigma (1), \sigma (2), \sigma (3)} 
\subset SL(6)$ 
of flag types $(\sigma (1), \sigma (2), 
\sigma (3))$ where $\sigma$ are permutations 
of $\{1,2,3\}$. Put $Y_{i,j,k} := 
T^*(SL(6)/P_{i,j,k})$. 
They are symplectic resolutions of 
$\bar{\mathcal{O}}$.  
Then we have 
a link of birational maps: 
$$  --\to Y_{321} --\to Y_{231} --\to Y_{213} --\to $$ 
$$ Y_{123} --\to Y_{132} --\to Y_{312} --\to Y_{321}.$$ 

Each birational map fits into a diagram which is locally 
a trivial family of Mukai flops of type $A$. 
For example, $Y_{321}$ and $Y_{231}$ are linked by 
a diagram 
$$ Y_{321} \rightarrow X_{5,1} \leftarrow Y_{231}. $$ 
This diagram is locally a trivial family of the 
Mukai flop 
$$ T^*G(2,5) \rightarrow \bar{\mathcal{O}}_{[2^2,1]} 
\leftarrow T^*G(3,5). $$
Let $N^1(Y_{i,j,k})$ be the 
Abelian groups of $\bar{\mathcal{O}}$-numerical classes of 
$\mathbf{R}$-divisors on $Y_{i,j,k}$. Note that 
$N^1(Y_{i,j,k}) \cong \mathbf{R}^2$. 
Since $Y_{i,j,k}$ are isomorphic in codimension 1, 
$N^1(Y_{i,j,k})$ are naturally identified. Hence we 
denote by $N^1$ these $\mathbf{R}$-vector spaces. 
$N^1$ is divided into six chambers, each of which 
corresponds to the ample cone 
$\overline{\mathrm{Amp}}(Y_{i,j,k})$.    
\end{Exam} 

\begin{Exam}
Let $\mathcal{O} \subset \mathfrak{so}(10)$ be 
the nilpotent orbit of Jordan type 
$[4^2, 1^2]$. Take an element $x$ of $\mathcal{O}$. 
Then $x$ has four polarizations 
$P^+_{3,2,2,3}$, $P^-_{3,2,2,3}$, $P^+_{2,3,3,2}$ 
and $P^-_{2,3,3,2}$, each of which has the flag 
type indicated by the index. Note that there are 
different polarizations which have the same 
flag type. We put $Y^+_{i,j,j,i} := 
T^*(SO(10)/P^+_{i,j,j,i})$ and        
$Y^-_{i,j,j,i} := 
T^*(SO(10)/P^-_{i,j,j,i})$.  
They are symplectic resolutions of $\bar{\mathcal{O}}$. 
Then we have a link of birational maps 
$$ Y^+_{3,2,2,3} --\to Y^+_{2,3,3,2} --\to 
Y^-_{2,3,3,2} --\to Y^-_{3,2,2,3}. $$ 
Here $Y^+_{3,2,2,3}$ and $Y^+_{2,3,3,2}$ are linked by 
a diagram 
$$ Y^+_{3,2,2,3} \rightarrow X^+_{5,5} 
\leftarrow Y^+_{2,3,3,2}.$$ 
This diagram is locally a trivial family of 
the Mukai flop 
$$ T^*G(3,5) \rightarrow \bar{\mathcal{O}}_{[2^2,1]} 
\leftarrow T^*G(2,5). $$
The same picture can be drawn for 
$Y^-_{3,2,2,3}$ and $Y^-_{2,3,3,2}$. 
Finally, $Y^+_{2,3,3,2}$ and $Y^-_{2,3,3,2}$ 
are linked by a diagram 
$$ Y^+_{2,3,3,2} \rightarrow X^+_{2,6,2} 
\leftarrow Y^-_{2,3,3,2}.$$ 
This diagram is locally a trivial family of 
the Mukai flop 
$$ T^*G^+_{iso}(3,6) 
\rightarrow \bar{\mathcal{O}}_{[2^2,1^2]} 
\leftarrow T^*G^-_{iso}(3,6). $$  
$Y^+_{2,3,3,2}$ (resp. $Y^-_{2,3,3,2}$) 
has two small birational contraction maps. 
On the other hand, $Y^+_{3,2,2,3}$ 
(resp. $Y^-_{3,2,2,3}$) has a small 
birational contraction map and a divisorial 
birational contraction map. This is the reason 
why the link of birational maps above ends 
at $Y^+_{3,2,2,3}$ and $Y^-_{3,2,2,3}$. 
Let $N^1(Y^+_{i,j,j,i})$ 
(resp. $N^1(Y^-_{i,j,j,i})$) be the 
$\mathbf{R}$-vector space of  
$\bar{\mathcal{O}}$-numerical classes of 
$\mathbf{R}$-divisors on $Y^+_{i,j,j,i}$ 
(resp. $Y^-_{i,j,j,i}$). Note that 
these spaces have dimension 2. 
Since $Y^+_{i,j,j,i}$ and $Y^-_{i,j,j,i}$ 
are isomorphic in codimension 1, 
these $\mathbf{R}$-vector spaces are naturally 
identified. Hence we 
denote them by $N^1$. 
Let $\mathrm{Mov} \subset N^1$ be the 
subcone generated by movable $\mathbf{R}$ 
divisors. In our case $\mathrm{Mov}$ does not 
coincides with $N^1$. 
The closure $\overline{\mathrm{Mov}}$ of 
$\mathrm{Mov}$ is divided into four chambers, 
each of which 
corresponds to the ample cone 
$\overline{\mathrm{Amp}}(Y^+_{i,j,j,i})$ 
or $\overline{\mathrm{Amp}}(Y^-_{i,j,j,i})$.
\end{Exam}     

\begin{Exam}
Let $\mathcal{O} \subset \mathfrak{so}(10)$ be 
the nilpotent orbit of Jordan type 
$[3^2, 2^2]$. Take an element $x$ of $\mathcal{O}$. 
Then $x$ has three polarizations 
$P^+_{1,4,4,1}$, $P^-_{1,4,4,1}$ 
and $P_{4,1,1,4}$, each of which has the flag 
type indicated by the index. Note that there are 
different polarizations which have the same 
flag type. We put $Y^+_{1,4,4,1} := 
T^*(SO(10)/P^+_{1,4,4,1})$,         
$Y^-_{1,4,4,1} := 
T^*(SO(10)/P^-_{1,4,4,1})$ 
and $Y_{4,1,1,4} := T^*(SO(10)/P_{4,1,1,4})$.   
They are symplectic resolutions of $\bar{\mathcal{O}}$. 
Then we have a link of birational maps 
$$ Y^+_{1,4,4,1} --\to Y_{4,1,1,4} --\to 
Y^-_{1,4,4,1}. $$ 
Here $Y^+_{1,4,4,1}$ and $Y_{4,1,1,4}$ are linked by 
a diagram 
$$ Y^+_{1,4,4,1} \rightarrow X^+_{5,5} 
\leftarrow Y_{4,1,1,4}.$$ 
This diagram is locally a trivial family of 
the Mukai flop 
$$ T^*G(1,5) \rightarrow \bar{\mathcal{O}}_{[2,1^3]} 
\leftarrow T^*G(4,5). $$
The same picture can be drawn for 
$Y^-_{1,4,4,1}$ and $Y_{4,1,1,4}$.   
$Y_{4,1,1,4}$  
has two small birational contraction maps. 
On the other hand, $Y^+_{1,4,4,1}$ 
(resp. $Y^-_{1,4,4,1}$) has a small 
birational contraction map and a divisorial 
birational contraction map. This is the reason 
why the link of birational maps above ends 
at $Y^+_{1,4,4,1}$ and $Y^-_{1,4,4,1}$. 
Let $N^1(Y^+_{1,4,4,1})$   
(resp. $N^1(Y^-_{1,4,4,1})$,  
$N^1(Y_{4,1,1,4})$ ) be the $\mathbf{R}$-vector space 
of $\bar{\mathcal{O}}$-numerical classes of $\mathbf{R}$   
divisors of $Y^+_{1,4,4,1}$ (resp. $Y^-_{1,4,4,1}$,  
$Y_{4,1,1,4}$).   
Since $Y^+_{1,4,4,1}$, $Y^-_{1,4,4,1}$ and 
$Y_{4,1,1,4}$ 
are isomorphic in codimension 1, 
these $\mathbf{R}$-vector spaces are naturally 
identified. Hence we 
denote them by $N^1$. 
Let $\mathrm{Mov} \subset N^1$ be the 
subcone generated by movable $\mathbf{R}$ 
divisors. In our case $\mathrm{Mov}$ does not 
coincides with $N^1$. 
The closure $\overline{\mathrm{Mov}}$ of 
$\mathrm{Mov}$ is divided into three chambers. 
They correspond to the ample cones of 
three resolutions of $\bar{\mathcal{O}}$.   
\end{Exam}

The following is a special case of a more general 
problem (cf. \cite{Na 2}).  

\begin{Conj}
All symplectic resolutions of a nilpotent orbit 
closure in a classical Lie algebra have equivalent 
bounded derived categories of coherent sheaves. 
\end{Conj} 

In the conjecture, the equivalences between derived 
categories should respect the birational 
map between resolutions. More explicitly, let 
$Y$ and $Y'$ be two symplectic resolutions and let 
$U \subset Y$ and $U' \subset Y'$ be Zariski open 
subsets such that $U$ is isomorphically mapped onto $U'$ 
by the natural birational map $\phi : Y \to Y'$. Then 
we want to have an equivalece 
$F: D(Y) \to D(Y')$ such that 
$F(\mathcal{O}_y) \cong \mathcal{O}_{\phi(y)}$ for 
all $y \in U$, where $\mathcal{O}_y$ (resp. $\mathcal{O}
_{\phi(y)}$) is the structure sheaf of the closed point 
$y$ (resp. $\phi(y)$).  
By Theorem 4.4 and by \cite{Na 1}, \S 5, the conjecture 
is reduced to the cases of Mukai flops of type $A$ and 
of type $D$.     
\vspace{0.12cm}

\section{Deformations of nilpotent orbits}
Let $x \in \g$ be a nilpotent element of 
a Lie algebra attached to a classical simple 
complex Lie group $G$.  
Let $\mathcal{O}$ be the nilpotent orbit 
containing $x$. In this section, by using 
an idea of Borho and Kraft \cite{B-K}, we shall 
construct a {\em deformation}\footnote
{Here we do not assume that $f$ is flat, but 
only assume that all fibers of $f$ has the 
same dimension. This suffices for an application.} 
$f: \mathcal{S} 
\to \mathfrak{k}$ of $\bar{\mathcal{O}}$ 
such that 
\vspace{0.12cm}

(i) $f^{-1}(0) = \bar{\mathcal{O}}$ for 
$0 \in \mathfrak{k}$, and 
\vspace{0.12cm}

(ii) for any Springer resolution $T^*(G/P) \to 
\bar{\mathcal{O}}$, there is a {\em simultaneous 
resolution} 
$$ \tau_P: \mathcal{E}_P \to \mathfrak{k}$$ 
of $f$, where $(\tau_P)^{-1}(0) = 
T^*(G/P)$ and where $(\tau_P)^{-1}(t) 
\cong f^{-1}(t)$ for a general point 
$t \in \mathfrak{k}$. 
\vspace{0.12cm}

As a corollary of this construction, we 
can verify Conjecture 2 in \cite{F-N} for the 
closure of a nilpotent orbit of a classical 
simple Lie algebra. 
Conjecture 2 has already been proved for 
$\mathfrak{sl}(n)$ in \cite{F-N}, Theorem 4.4 
in a very explicit form. 
Note that, a weaker version 
of this conjecture has been proved by Fu \cite{Fu 2} 
for the closure of a nilpotent orbit of a 
classical simple Lie algebra. 

$\g$ becomes a $G$-variety via the adjoint action. 
Let $Z \subset \g$ be a closed subvariety. 
For $m \in \mathbf{N}$, put 
$$ Z^{(m)} := \{x \in Z; \dim Gx = m\}. $$ 
$Z^{(m)}$ becomes a locally closed subset of 
$Z$. We put $m(Z) := \mathrm{max}\{m; 
m = \dim Gx, \exists x \in Z\}$. Then 
$Z^{m(Z)}$ is an open subset of $Z$, which will 
be denoted by $Z^{\mathrm{reg}}$. 
A {\em sheet} of $Z$ is an irreducible 
component of some $Z^{(m)}$.      
A sheet of $\g$ is called a {\em Dixmier sheet} 
if it contains a semi-simple element of $\g$. 

Let $P \subset G$ be a parabolic subgroup and 
let $\mathfrak{p}$ be its Lie algebra. 
Let $\mathfrak{m}(P)$ be 
the Levi factor of $\mathfrak{p}$. 
We put $\mathfrak{k}(P) := \g^{\mathfrak{m}(P)}$ 
where 
$$ \g^{\mathfrak{m}(P)} := \{x \in \g; 
[x,y] = 0, \forall y \in \mathfrak{m}(P)\}. $$ 
Let $\mathfrak{r}(P)$ be the radical of $\mathfrak{p}$. 

\begin{Thm}
$G\mathfrak{r}(P) = \overline{G\mathfrak{k}(P)}$ 
and $G\mathfrak{r}(P)^{\mathrm{reg}}$ $(= 
\overline{G\mathfrak{k}(P)}^{\mathrm{reg}})$ is 
a Dixmier sheet. 
\end{Thm}   
   
{\em Proof}. See \cite{B-K}, Satz 5.6. 
\vspace{0.12cm}

Every element $x$ of $\g$ can be uniquely written as 
$x = x_n + x_s$ with $x_n$ nilpotent and 
with $x_s$ semi-simple 
such that $[x_n,x_s] = 0$. 
Let $\mathfrak{h}$ be a Cartan subalgebra and 
let $W$ be the Weyl group with respect to 
$\mathfrak{h}$. The set of semi-simple orbits is 
identified with $\mathfrak{h}/W$. Let 
$\g \to \mathfrak{h}/W$ be the map defined as 
$x \to [\mathcal{O}_{x_s}]$. 
There is 
a direct sum decomposition 
$$ \mathfrak{r}(P) = \mathfrak{k}(P) 
\oplus \mathfrak{n}(P), (x \to x_1 + x_2)$$ 
where $\mathfrak{n}(P)$ is the 
nil-radical of $\mathfrak{p}$ (cf. \cite{Slo}, 
4.3). 
We have a well-defined map 
$$ G \times_P \mathfrak{r}(P) \to 
\mathfrak{k}(P) $$ 
by sending 
$[g, x] \in G \times_P \mathfrak{r}(P)$ to 
$x_1 \in \mathfrak{k}(P)$ and  
there is a commutative diagram 
$$ G \times_P \mathfrak{r}(P) \to G\mathfrak{r}(P) $$ 
$$ \downarrow \hspace{1.0cm} \downarrow $$ 
$$ \mathfrak{k}(P) \to \mathfrak{h}/W. $$ by 
\cite{Slo}, 4.3.   
           
\begin{Lem}
The induced map 
$$ G \times_P \mathfrak{r}(P) \stackrel{\mu_P}\to 
\mathfrak{k}(P) \times_{\mathfrak{h}/W} 
G\mathfrak{r}(P)$$ 
is a birational map. 
\end{Lem} 

{\em Proof}. Let $h \in 
\mathfrak{k}(P)^{\mathrm{reg}}$ and 
denote by $\bar{h} \in \mathfrak{h}/W$ its 
image by the map $\mathfrak{k}(P) \to 
\mathfrak{h}/W$. Then the fiber of 
the map $G\mathfrak{r}(P) \to 
\mathfrak{h}/W$ over $\bar{h}$ coincides 
with a semi-simple orbit $\mathcal{O}_h$ of $\g$ 
containing $h$. In fact, by  
Theorem 5.1, the fiber actally contains 
this orbit. The fiber is closed in $\g$ 
because $G\mathfrak{r}(P)$ is closed 
subset of $\g$ by Theorem 5.1. Note 
that a semi-simple orbit of $\g$ is 
also closed. Hence if the fiber and 
$\mathcal{O}_h$ does not coincide, 
then the fiber contains an orbit with 
larger dimension than $\dim \mathcal{O}_h$. 
This contradicts the fact that 
$\overline{G\mathfrak{k}(P)^{\mathrm{reg}}} 
= G\mathfrak{r}(P)$.  
Take a point 
$(h,h') \in \mathfrak{k}(P)^{\mathrm{reg}} 
\times_{\mathfrak{h}/W} G\mathfrak{r}(P)$. 
Then $h'$ is a semi-simple element 
$G$-conjugate to $h$. 
Fix an element $g_0 \in G$ such that 
$h' = g_0h(g_0)^{-1}$. 
We have 
$$ (\mu_P)^{-1}(h,h') 
= \{[g,x] \in G \times_P \mathfrak{r}(P);\, 
x_1 = h,\, gxg^{-1} = h'\}. $$ 
Since $x = px_1p^{-1}$ for some $p \in P$ and 
conversely $(px_1p^{-1})_1 = x_1$ for any $p \in P$ 
(cf. \cite{Slo}, Lemma 2, p.48), we have 
$$ (\mu_P)^{-1}(h,h') = 
\{[g,php^{-1}] \in G \times_P\mathfrak{r}(P);\, 
g \in G,\, p \in P,\, (gp)h(gp)^{-1} = h'\} =$$ 
$$ \{[gp,h] \in G \times_P\mathfrak{r}(P);\, 
g \in G,\, p \in P,\, (gp)h(gp)^{-1} 
= h'\} = $$ 
$$ \{[g,h] \in G \times_P\mathfrak{r}(P);\, 
ghg^{-1} = h' \} = $$  
$$ \{[g_0g',h] \in G \times_P \mathfrak{r}(P);\, 
g' \in Z_G(h)\} = g_0(Z_G(h)/Z_P(h)). $$ 
Here $Z_G(h)$ (resp. $Z_P(h)$) is the 
centralizer of $h$ in $G$ (resp. $P$). 
By \cite{Ko}, 3.2, Lemma 5, $Z_G(h)$  
is connected. Moreover, since 
$\g^h \subset \mathfrak{p}$, 
$\mathrm{Lie}(Z_G(h)) = 
\mathrm{Lie}(Z_P(h))$. 
Therefore, $Z_G(h)/Z_P(h) 
= \{1\}$, and $(\mu_P)^{-1}(h,h')$ 
consists of one point. 

\begin{Lem}
The map 
$$ G \times_P \mathfrak{r}(P) \to 
G\mathfrak{r}(P) $$ 
is a proper map.
\end{Lem}

{\em Proof}. 
As a vector subbundle, we have a closed 
immersion 
$$G \times_P \mathfrak{r}(P) \to 
G/P \times \g. $$ 
This map factors through 
$G/P \times G\mathfrak{r}(P)$, and hence 
we have a closed immersion 
$$G \times_P \mathfrak{r}(P) \to 
G/P \times G\mathfrak{r}(P).$$ 
Our map is the composition of 
this closed immersion and the 
projection 
$$ G/P \times G\mathfrak{r}(P) 
\to G\mathfrak{r}(P).$$  Since 
$G/P$ is compact, this projection 
is a proper map.  

\begin{Lem}
Let $x \in \g$ be a nilpotent 
orbit of a classical simple Lie 
algebra and denote by $\mathcal{O}$ 
the nilpotent orbit containing 
$x$. Then the polarizations of 
$x$ giving Springer resolutions 
of $\bar{\mathcal{O}}$ all have conjugate 
Levi factors. 
\end{Lem}

{\em Proof}. By Theorem 2.8, the Levi 
type $\pi$ of such a polarization $P$ 
is unique. 
Parabolic 
subgroups of $G$ with the same Levi 
type have congugate Levi factors except in 
the case $q = \epsilon = 0$ and 
$\pi^i \equiv 0$ (mod 2) $\forall i$ 
(cf. \cite{He}, Lemma 4.6,(c)). 
In this exceptional case, 
there are two conjugacy classes of parabolic subgroups 
having non-conjugate 
Levi factors. They have mutually 
different Richardson orbits 
$\mathcal{O}^I$ and $\mathcal{O}^{II}$. 
Here $\mathcal{O}^I$ and 
$\mathcal{O}^{II}$ are very even orbits 
with the same Jordan type. In particular, 
it is impossible that two parabolic subgroups 
having non-conjugate Levi factors becomes 
polarizations of the same element $x$.   

\begin{Lem}
Let $x \in \g$ be the same as the 
previous lemma. Let  
$P$ and $P'$ be polarizations of $x$. 
Assume that they both give Springer resolutions 
of $\bar{\mathcal{O}}$. 
Then $\mathfrak{k}(P)$ and 
$\mathfrak{k}(P')$ are conjugate to 
each other. 
\end{Lem}

{\em Proof}. 
Let $M_P$ and $M_{P'}$ be 
Levi factors of $P$ and $P'$ 
respectively. Then 
$M_P$ and $M_{P'}$ are conjugate by 
the previous lemma. Hence their 
centralizers are also conjugate. 
The Lie algebras of these centralizers 
are $\mathfrak{k}(P)$ and 
$\mathfrak{k}(P')$.  

\begin{Cor}
For $P$, $P' \in \mathrm{Pol}(x)$ which 
give Springer resolutions of 
$\bar{\mathcal{O}}$, we have $G\mathfrak{r}(P) 
= G\mathfrak{r}(P')$.   
\end{Cor}

{\em Proof}. 
By Theorem 5.1, $G\mathfrak{r}(P) 
= \overline{G\mathfrak{k}(P)}$ and    
$G\mathfrak{r}(P') 
= \overline{G\mathfrak{k}(P')}$. 
Since $G\mathfrak{k}(P) = 
G\mathfrak{k}(P')$, we have the 
result. 

\begin{Lem}
The image of the composed map 
$$G\mathfrak{r}(P) \to \g \to 
\mathfrak{h}/W$$ 
coincides with 
$\mathfrak{k}(P)/W_P$, where 
$$ W_P = \{w \in W; w(\mathfrak{k}(P)) 
= \mathfrak{k}(P)\}. $$ 
\end{Lem}

{\em Proof}. By definition, 
$\mathfrak{k}(P)/W_P \subset 
\mathfrak{h}/W$, which is a 
closed subset. Since 
$G\mathfrak{r}(P) = 
\overline{G\mathfrak{k}(P)}$, 
we only have to prove that 
the image of $G\mathfrak{k}(P)$ 
by the map $\g \to \mathfrak{h}/W$ 
coincides with $\mathfrak{k}(P)/W_P$. 
Every element of $G\mathfrak{k}(P)$ 
is semi-simple, and the map 
$G\mathfrak{k}(P) \to \mathfrak{h}/W$ 
sends an element of $G\mathfrak{k}(P)$ 
to its (semi-simple) orbit. 
Hence the image coincides with 
$\mathfrak{k}(P)/W_P$.  

\begin{Cor}
$\mathfrak{k}(P)$ and $\mathfrak{k}(P')$ 
are $W$-conjugate in $\mathfrak{h}$.  
\end{Cor}

{\em Proof}.    
Let $q : \mathfrak{h} \to \mathfrak{h}/W$ 
be the quotient map. 
Since $G\mathfrak{r}(P) = 
G\mathfrak{r}(P')$, $q(\mathfrak{k}(P)) 
= q(\mathfrak{k}(P'))$ by the previous 
lemma. Put $\mathfrak{k}_{\pi} := 
q(\mathfrak{k}(P))$. Then $\mathfrak{k}(P)$ 
and $\mathfrak{k}(P')$ are both irreducible 
components of $q^{-1}(\mathfrak{k}_{\pi})$. 
Hence, $\mathfrak{k}(P)$ and $\mathfrak{k}(P')$ 
are $W$-conjugate in $\mathfrak{h}$. 
\vspace{0.12cm}

We fix a polarization $P_0$ of $x$ 
which gives a Springer resolution of 
$\bar{\mathcal{O}}$. Let $P$ be another 
such polarization. By the corollary above, 
$\mathfrak{k}(P)$ and $\mathfrak{k}(P_0)$ 
are finite coverings of $\mathfrak{k}_{\pi}$, and  
there is a $\mathfrak{k}_{\pi}$-isomorphism 
$\mathfrak{k}(P) \cong \mathfrak{k}(P_0)$. 
We fix such an isomorphism. Then 
it induces an isomorphism 
$$ \mathfrak{k}(P)\times_{\mathfrak{h}/W}G\mathfrak{r}(P) 
\stackrel{\iota_P}\to  
\mathfrak{k}(P_0)\times_{\mathfrak{h}/W}
G\mathfrak{r}(P_0). $$ 
We put $\nu_P := \iota_P \circ \mu_P$, and  
$$ \mathcal{S} := \mathfrak{k}(P_0)\times_{\mathfrak{h}/W}
G\mathfrak{r}(P_0). $$ 
Denote by $f$ the first projection  
$\mathcal{S} \to \mathfrak{k}(P_0)$.  
Then $f^{-1}(0) = \bar{\mathcal{O}}$, and 
for each polarization $P$ of $x$, 
$$ G \times_P \mathfrak{r}(P) 
\stackrel{\nu_P}\to \mathcal{S} \to 
\mathfrak{k}(P_0)$$ 
gives a simultaneous resolution of $f$.  
This simultaneous resolution coincides 
with the Springer resolution 
$T^*(G/P) \to \bar{\mathcal{O}}$ over 
$0 \in \mathfrak{k}(P_0)$.          
\vspace{0.12cm}

The following conjecture is posed in 
\cite{F-N}. 

\begin{Conj}
Let $Z$ be a normal symplectic singularity. 
Then for any two symplectic resolutions 
$f_i: X_i \to Z$, $i=1,2$, there are flat 
deformations $\mathcal{X}_i 
\stackrel{F_i}\to \mathcal{Z} 
\to T$ such that, for $t \in T - \{0\}$, 
$F_{i,t}: \mathcal{X}_{i,t} \to 
\mathcal{Z}_t$ are isomorphisms. 
\end{Conj} 

\begin{Thm}
The conjecture holds for the normalization 
$\tilde{\mathcal{O}}$ of a nilpotent orbit closure 
$\bar{\mathcal{O}}$ in a classical Lie algebra.    
\end{Thm}

{\em Proof}. By \cite{Fu 1}, Theorem 3.3, all symplectic 
resolutions of $\tilde{\mathcal{O}}$ are Springer 
resolutions. Take a general curve $T \subset 
\mathfrak{k}(P_0)$ passing through $0 \in 
\mathfrak{k}(P_0)$, and pull back the family 
$$G \times_P \mathfrak{r}(P) \stackrel{\nu_P}\to 
\mathcal{S} \to \mathfrak{k}(P_0)$$ 
by $T \to \mathfrak{k}(P_0)$. Put $\bar{\mathcal{Z}} 
:= \mathcal{S} \times_{\mathfrak{k}(P_0)}T$. Then, 
for each $P$, we have a simultaneous resolution of 
$\bar{\mathcal{Z}} \to T$: 
$$ \mathcal{X}_P \to \bar{\mathcal{Z}} \to T. $$ 
Let $\mathcal{Z}$ be the normalization of 
$\bar{\mathcal{Z}}$. 
Then the map $\mathcal{X}_P \to \bar{\mathcal{Z}}$ 
factors through $\mathcal{Z}$. Now 
$$ \mathcal{X}_P \to \mathcal{Z} \to T$$ 
gives a desired deformation of the Springer 
resolution $T^*(G/P) \to \tilde{\mathcal{O}}$.       

\begin{Exam}
Our abstract construction coincides with the 
explicit construction in \cite{F-N}, Theorem 4.4 
in the case where $\g = \mathfrak{sl}(n)$.  
Let us briefly observe the correspondence between 
two constructions. Assume that $\mathcal{O}_x \subset 
\mathfrak{sl}(n)$ is the orbit containing an nilpotent 
element $x$ of type 
$\mathbf{d} := [d_1, ..., d_k]$. 
Let $[s_1, ..., s_m]$ be the 
dual partition of $\mathbf{d}$ (cf. Notation and 
Convention). By Theorem 2.6, the polarizations of 
$x$ have the flag type $(s_{\sigma (1)}, \cdots, 
s_{\sigma (m)})$ with $\sigma \in \Sigma_m$. 
We denote them by $P_{\sigma}$.  
We put $P_0 := P_{id}$.  
Define $F_\sigma : = 
SL(n)/P_{\sigma}. $ 
Let $$ \tau_1  \subset \cdots \subset  \tau_{m-1} \subset
\cit^n \otimes_{\cit}\mathcal{O}_{F_{\sigma}}$$ 
be the universal subbundles on $F_{\sigma}$. 
A point of $T^*F_\sigma$ is expressed as a pair 
$(p, \phi)$ of $p \in F_{\sigma}$ and $\phi \in 
\mathrm{End}(\cit^n)$ such that 
$$ \phi(\cit^n) \subset \tau_{m-1}(p), \cdots,  
\phi(\tau_2(p)) \subset 
\tau_1(p), \phi(\tau_1(p)) = 0.$$ 
The Springer resolution 
$$ s_{\sigma}: T^*F_{\sigma} \to \bar{\mathcal{O}}$$ 
is defined as $s_{\sigma}((p, \phi)) := \phi$.
In \cite{F-N}, Theorem 4.4, we have next defined 
a vector bundle $\mathcal{E}_{\sigma}$ with an 
exact sequence 
$$ 0 \to T^*F_{\sigma} \to \mathcal{E}_{\sigma} 
\stackrel{\eta_{\sigma}}\to 
\mathcal{O}_{F_{\sigma}}^{m-1} \to 0.$$ 
For $p \in F_{\sigma}$, we can choose a basis 
of $\mathbf{C}^n$ such that $T^*F_{\sigma}(p)$ 
consists of the matrices of the following form 
 
$$ 
\begin{pmatrix} 
0 & * & \cdots &* \\
0 & 0 &  \cdots & *\\
\cdots & & & \cdots \\
0 & 0 &  \cdots& 0 \\ 
\end{pmatrix}.
$$  

Then $\mathcal{E}_{\sigma}(p)$ is the vector 
subspace of $\mathfrak{sl}(n)$ consisting of 
the matrices $A$ of the following 
form
$$ 
\begin{pmatrix} 
a_{\sigma (1)} & * & \cdots & * \\
0 & a_{\sigma (2)} &\cdots  & *\\
\cdots & & & \cdots \\
0 & 0 & \cdots  &a_{\sigma (m)} \\ 
\end{pmatrix},$$
where $a_i := a_i I_{s_i}$ and $I_{s_i}$ is 
the identity matrix of 
the size $s_i \times s_i$. 
Since $A \in \mathfrak{sl}(n)$, 
$\Sigma_i s_i a_i = 0$. Here we define the 
map $\eta_{\sigma}(p): \mathcal{E}_{\sigma}(p) 
\to {\cit}^{\oplus m-1}$ as 
$\eta_{\sigma}(p)(A) := (a_1, a_2, \cdots, a_{m-1}).$  
This vector bundle $\mathcal{E}_{\sigma}$ 
is nothing but our 
$SL(n) \times_{P_{\sigma}}
\mathfrak{r}(P_{\sigma})$.   
Moreover, the map 
$$\eta_{\sigma}: 
\mathcal{E}_{\sigma} \to \cit^{m-1} $$ 
coincides with the map 
$$ SL(n) \times_{P_{\sigma}}
\mathfrak{r}(P_{\sigma}) \to 
\mathfrak{k}(P_0), $$ 
where we identify $\mathfrak{k}(P_{\sigma})$ 
with $\mathfrak{k}(P_0)$ by an 
$\mathfrak{k}_{\pi}$-isomorphism.   
Finally, in \cite{F-N}, Theorem 4.4 we have 
defined $\overline{N} \subset \mathfrak{sl}(n)$ 
to be the set of 
all matrices which is conjugate to a matrix 
of the following form:
$$ 
\begin{pmatrix} 
b_{1} & * & \cdots & * \\
0 & b_{2} &\cdots  & *\\
\cdots & & & \cdots \\
0 & 0 & \cdots  & b_{m} \\ 
\end{pmatrix},
$$
where $b_i = b_i I_{s_i}$ and  
$I_{s_i}$ is the identity matrix of order $s_{i}$.  
Furthermore 
the zero trace 
condition $\sum_i s_i b_i =0$ was required.   
For $A \in \overline{N}$, let $\phi_A(x) := \mathrm{det}
(xI - A)$ be the characteristic polynomial of 
$A$. Let $\phi_i(A)$ be the coefficient of 
$x^{n-i}$ in $\phi (A)$. 
Here the characteristic map  
$ch: \overline{N} 
\to \cit^{n-1}$ has been defined as 
$ch(A) := (\phi_2(A), ..., 
\phi_n(A))$. 
This $\bar{N}$ is nothing but our 
$SL(n)\mathfrak{r}(P_{\sigma})$. 
As is proved in Corollary 5.6, this is 
independent of the choice of $P_{\sigma}$. 
The characteristic map $ch$ above coincides 
with the composed map 
$$ SL(n)\mathfrak{r}(P_{\sigma}) 
\subset \mathfrak{sl}(n) \to \mathfrak{h}/W. $$
\end{Exam}

\quad \\
\quad\\

Yoshinori Namikawa \\
Departement of Mathematics, 
Graduate School of Science, Kyoto University,
Kita-shirakawa Oiwake-cho, Kyoto, 606-8502, JAPAN \\
namikawa@kusm.kyoto-u.ac.jp

 　　　　      
\end{document}